# SPGD: Steepest Perturbed Gradient Descent Optimization

Amir M. Vahedi and Horea T. Ilies
School of Mechanical, Aerospace, and Manufacturing Engineering
University of Connecticut

December 28, 2025

**Abstract**

Optimization algorithms are pivotal in advancing various scientific and industrial fields but often encounter obstacles such as trapping in local minima, saddle points, and plateaus (flat regions), which makes the convergence to reasonable or near-optimal solutions particularly challenging.

This paper presents the Steepest Perturbed Gradient Descent (SPGD), a novel algorithm that innovatively combines the principles of the gradient descent method with periodic uniform perturbation sampling to effectively circumvent these impediments and lead to better solutions whenever possible. SPGD is distinctively designed to generate a set of candidate solutions and select the one exhibiting the steepest loss difference relative to the current solution. It enhances the traditional gradient descent approach by integrating a strategic exploration mechanism that significantly increases the likelihood of escaping sub-optimal local minima and navigating complex optimization landscapes effectively. Our approach not only retains the directed efficiency of gradient descent but also leverages the exploratory benefits of stochastic perturbations, thus enabling a more comprehensive search for global optima across diverse problem spaces. We demonstrate the efficacy of SPGD in solving the 3D component packing problem, an NP-hard challenge. Preliminary results show a substantial improvement over six established methods, particularly on response surfaces with complex topographies and in multidimensional non-convex continuous optimization problems. Comparative analyses with established 2D benchmark functions over 30 randomized initial points highlight SPGD's robustness and reliability in non-convex optimization. These results emphasize SPGD's potential as a versatile tool for a wide range of optimization problems.

## 1 Introduction

Mathematical optimization is a fundamental process in engineering, science, and economics. Its main objective is to find solutions that minimize a predefined objective, typically expressed in terms of a real-valued function, while adhering to given constraints. This pursuit of optimal solutions is crucial in solving complex problems, where achieving the best possible results necessitates a careful balance of numerous factors and variables.

Among the many optimization techniques available, the gradient descent (GD) method stands out as a foundational and extensively used tool, and its origins can be traced back to Cauchy's pioneering work [1]. However, despite its widespread use, the gradient descent method has certain limitations. One of its major drawbacks is its tendency to get trapped in sub-optimal states, including saddle points and local minima, which may offer minimal improvement in solution quality. Additionally, the method may encounter difficulties in making progress towards the desired outcome when faced with flat regions in the problem space.

To address these challenges, extensive research efforts have been focused on enhancing the performance of the gradient descent method. As a result, numerous variants have been developed, each specifically designed to overcome the aforementioned pitfalls [2]. One notable variant is the Perturbed Gradient Descent (PGD), which has gained attention for its ability to navigate away from saddle points and potentially converge towards second-order optimal points [3].

In this paper, we present a strategic randomized perturbation algorithm combined with the gradient descent method, leveraging the strengths of both exploring the search space through randomized perturbation and converging to optimal points using gradient information. By introducing cyclical perturbations, our approach strategically balances the need for exploration with the efficiency of exploitation. Moreover, applying perturbations periodically rather than at every iteration significantly reduces computational costs, making the optimization process more efficient without sacrificing the thoroughness of the search. It promises a more reliable pathway to discovering superior solutions, thereby expanding the horizon of possibilities in optimization challenges. This enhanced method is designed not only to navigate more effectively through the complexities of practical optimization landscapes but also to refine the search for optimal solutions with greater precision.

The remainder of this paper systematically explores the Steepest Perturbed Gradient Descent (SPGD) algorithm and its comparative advantages in the domain of optimization. Section 2 delves into a variety of related methodologies, focusing on variants of the gradient descent method and the integration of perturbation sampling techniques. These approaches establish a foundation for understanding the landscape of optimization strategies and highlight the necessity for innovations, such as SPGD. Section 3 is dedicated to a detailed exposition of the SPGD algorithm itself, including its theoretical underpinnings, algorithmic structure, and the rationale behind its design choices. Following this, Section 5 presents numerical results from a series of experiments designed to evaluate the performance of SPGD against various established optimization algorithms. These experiments are conducted on a selection of well-known optimization test functions, providing a rigorous comparison and demonstrating the practical implications of SPGD in addressing complex optimization challenges. Finally, Section 6 discusses the outcomes of these comparisons, emphasizing the superior performance of SPGD over the methods analyzed. The conclusions not only underscore the effectiveness of SPGD but also set the stage for future research directions and potential applications in broader optimization contexts.

## 2 Related Work

The gradient descent (GD) method is a first-order optimization algorithm that updates the design variables in the direction opposite to the gradient of the objective function with respect to those variables [2]. It's widely used due to its simplicity and efficiency in convex problems. The gradient descent method converges to a local optimal solution with a mathematical guarantee. Gradient descent tends to exploit local information to improve the solution iteratively. However, it may not explore the search space effectively, potentially getting trapped in local minima or saddle points, particularly in non-convex optimization landscapes. GD is known to struggle with flat areas where the gradient is close to zero, leading to slow or no progress [3].

Nesterov's Accelerated Gradient (NAG) method enhances traditional gradient descent by incorporating a forward-looking step. This tweak allows the optimizer to anticipate future gradients, reducing oscillations and speeding up convergence, particularly in convex settings. NAG is highly effective in training deep neural networks due to its efficiency in navigating high-dimensional data spaces. However, its performance can vary in non-convex environments with complex landscapes [4, 5]. For a comprehensive overview of gradient descent and its variants, we refer readers to [2], which synthesizes developments across machine learning and optimization literature.

Simulated annealing (SA) [6] and genetic algorithm [7] are heuristic sampling-based optimization algorithms that use randomness and selection mechanisms inspired by natural processes to explore the solution space and select the best candidates for further iteration. These algorithms can be used for different types of optimization problems, such as continuous, discrete, non-convex, and multi-objective problems [8, 9]. These methods may require significant computational resources and careful tuning of parameters (e.g., temperature in simulated annealing or mutation rate in evolutionary algorithms) to balance exploration and exploitation effectively.

Bayesian optimization (BO) is one of the sampling-based global optimization methods that has gained popularity, particularly in machine learning, for solving expensive black-box optimization problems. BO methods approximate the objective function using a surrogate probabilistic model, typically a Gaussian process (GP), which models the underlying function based on observed sample points [10, 11]. These methods balance exploration and exploitation by combining prior beliefs with posterior updates after each observation. The acquisition function, derived from this surrogate model, guides the search by quantifying the expected improvement or uncertainty in unexplored regions. Although BO is highly sample-efficient and effective in finding global optima with a limited number of function evaluations, especially in low-dimensional problems, it can be computationally expensive due to the cost of updating and optimizing the acquisition function at each iteration. Moreover, the performance of BO degrades in high-dimensional or highly non-smooth optimization landscapes [12]. A broader discussion of such sampling-based approaches can be found in [13], a recent survey on non-smooth optimization methods, including gradient sampling and probabilistic techniques.

Stochastic gradient descent (SGD) is a well-known optimization algorithm widely utilized in the training step of neural networks due to its efficiency in handling large datasets. Its capability to introduce randomness through mini-batches helps in escaping local minima, making it particularly suitable for large-scale machine learning problems where mini-batch sampling provides natural stochasticity [14]. The utilization of SGD presents an alternative view by adjusting parameters with a randomly chosen subset of data rather than the complete dataset, thereby injecting noise into the gradient approximations. By contrast, our SPGD is designed for deterministic nonconvex optimization problems where full gradients are available, and the challenge lies in escaping saddle points and local minima rather than handling stochastic gradients.

In the exploration of hybrid optimization methods, a notable approach combines the exploratory strengths of Simulated Annealing (SA) with the precise, local search capabilities of Gradient Descent (GD). This method strategically employs SA to break free from local optima by conducting a thorough search for a more promising solution candidate, upon which GD resumes. While this synergy offers a dynamic pathway to escape local minima, it introduces a significant computational burden. Moreover, this method diverges from traditional GD in that it cannot rely on the norm of the gradient as a criterion for termination. This alteration results in a less stringent stop condition, potentially affecting the algorithm's efficiency and termination reliability [15].

Another hybrid technique is perturbed gradient descent (PGD) that addresses the challenge of stagnation—a state where the gradient becomes negligible, and no further progress seems attainable in optimizing the objective function. This method introduces a single perturbation to the current solution when progress halts, effectively nudging the search process out of stagnation before proceeding with GD. This approach demonstrates an ability to escape saddle points effectively [3, 5, 16]. However, its performance is notably diminished in flat regions of the search space, where such perturbations fail to provide a meaningful direction for improvement.

In addition to the deterministic and heuristic methods previously discussed, the random walk method offers a stochastic approach to optimization that is particularly advantageous in complex, non-convex landscapes. Random walks operate by making a sequence of moves, each determined randomly in terms of direction and step size. This method inherently avoids the common pitfalls of gradient-based approaches, such as becoming trapped in local minima, by facilitating an unbiased exploration of the solution space. This characteristic is critical when dealing with high-dimensional optimization problems where the landscape is riddled with numerous local optima and saddle points [17]. Despite their potential for encompassing space exploration, random walks are often criticized for their inefficiency and slow convergence, especially in large-scale problems. They require a large number of iterations to approach the vicinity of a global optimum, as their exploration process lacks directionality inherent to methods like gradient descent or even simulated annealing. To address these limitations, researchers have explored hybrid strategies that combine the exploratory strengths of random walks with more systematic search techniques to balance exploration with exploitation more effectively [18, 19].

# 3 Methodology: SPGD

Traditional gradient descent algorithms efficiently exploit local gradient information to improve solutions iteratively. To minimize a given function $f : \mathbb{R}^d \rightarrow \mathbb{R}$, the updating rule at each iteration is [3]:

$$\mathbf{x}_{i+1} = \mathbf{x}_i - \alpha \nabla f(\mathbf{x}_i) \tag{1}$$

where $i$ is the number of current iteration, $\alpha > 0$ is the step size, and $\nabla f$ is the gradient of $f$.

However, in non-convex high-dimensional problems, the gradient descent method can become trapped in local minima or saddle points, missing out on globally optimal configurations. To address this limitation, we propose a novel algorithm that combines gradient descent with periodic randomized perturbations. These perturbations are particularly effective in non-convex, high-dimensional problems, where even small modifications can significantly alter the solution's position within the search space. This sensitivity to perturbations is crucial in navigating the complex terrain of such problems, where the landscape of potential solutions is riddled with local optima. By introducing strategically randomized perturbations, our algorithm enhances its ability to escape these local optima, thereby facilitating a more extensive exploration of the solution space. This periodic application of perturbations is key to avoiding the oscillatory behavior often observed in optimization trajectories of sampling-based methods, which can lead to inefficiencies and slow convergence. This approach becomes particularly advantageous in complex optimization scenarios characterized by challenges such as flatness, ruggedness, or saddle points of the objective surface, where conventional optimization algorithms might falter in making meaningful progress.

These perturbations are drawn from uniform random distributions[1] with constant amplitude profiles[2]. The uniform random distribution will create $N_P$ perturbed candidates around the gradient descent solution every $Iter_P$ iterations. All perturbed candidates will be evaluated and compared with the gradient descent solution. If the minimum value of perturbed candidates is equal or less than the value of the gradient descent solution, the corresponding perturbed candidate will be selected as the new solution. This policy of accepting solutions with equal objective values intentionally increases the algorithm's emphasis on exploration over exploitation within the optimization process. Such an approach is particularly advantageous in scenarios where the objective surface is flat, and traditional gradient descent methods stall due to insufficient gradient information. By facilitating exploration in these flat regions, SPGD ensures continued progress towards finding a global optimum, preventing the algorithm from becoming prematurely anchored to suboptimal solutions. The pseudo code of SPGD algorithm is described in Algorithm 1.

In the SPGD algorithm, the method of applying perturbations is adaptable to the specific requirements of the optimization problem at hand. For unconstrained optimization problems, such as 2D test function benchmarks and neural network training, perturbations are applied simultaneously to all variables, utilizing a uniform distribution with a constant amplitude. This ensures a broad, uniform exploration of the solution space, which is generally suitable for the landscapes presented by these types of problems.

However, when dealing with constrained problems like the 3D component packing, which present a complex optimization landscape, a different approach is warranted. In such scenarios, perturbation of a single variable can potentially lead to an infeasible solution or a worse candidate due to the constraints involved. To mitigate this, each variable is perturbed separately, effectively reducing the complexity and dimensionality of the optimization problem by focusing on one variable at a time, with all others held constant [20]. This targeted perturbation allows for a more controlled exploration of the solution space, ensuring that the search remains within feasible regions and is more likely to improve upon the current solution. This adaptive feature is designed to tailor the exploration process more precisely to the problem's landscape, enhancing the algorithm's flexibility and effectiveness in navigating constrained environments.

[1] Having a uniform distribution allows us to sample the solution space around the current solution uniformly within the range of $[-Amp, +Amp]$. This approach leads to more explorative sampling by equally covering the vicinity around the current position, rather than concentrating on the immediate area around the current solution or extending far beyond it.

[2] The main reason for choosing a constant amplitude profile is to maintain the simplicity of the algorithm in these benchmark functions. However, the amplitude profile ($Amp$) can be adjusted to any arbitrary profile based on the specific requirements of the optimization problem at hand.

**Algorithm 1** Steepest Perturbed Gradient Descent (SPGD)

**Input:** step size $\alpha > 0$, period $\mathrm{Iter}_P \in \mathbb{N}$, amplitude $\mathrm{Amp} > 0$, candidates $N_P \in \mathbb{N}$, horizon $\mathrm{Iter}_{\max}$
Initialize $\mathbf{x}_0$; set $i \leftarrow 0$
**while** $i < \mathrm{Iter}_{\max}$ **do**
  **if** $i \bmod \mathrm{Iter}_P \neq 0$ **then** ▷ plain GD
    $\mathbf{x}_{i+1} \leftarrow \mathbf{x}_i - \alpha \nabla f(\mathbf{x}_i)$
  **else** ▷ perturbation round (best-of-$N_P$)
    **for** $j = 1$ to $N_P$ **do**
      sample $\xi^{(j)} \sim \mathrm{Unif}(B_0(\mathrm{Amp}))$ ▷ ball of radius Amp
      $\mathbf{x}_i^{(j)} \leftarrow \mathbf{x}_i + \xi^{(j)}$
      $y_i^{(j)} \leftarrow f\left(\mathbf{x}_i^{(j)}\right)$
    **end for**
    $j^\star \in \arg\min_{j \in [N_P]} y_i^{(j)}$
    **if** $f(\mathbf{x}_i^{(j^\star)}) \leq f(\mathbf{x}_i)$ **then**
      $\mathbf{x}_{i+1} \leftarrow \mathbf{x}_i^{(j^\star)}$
    **end if**
  **end if**
  $i \leftarrow i + 1$
**end while**

The parameters of SPGD, notably the number of perturbations $N_P$, the perturbation interval $Iter_P$, and the perturbation amplitude $Amp$, play crucial roles in shaping the algorithm's behavior and performance. Increasing $N_P$ enhances the likelihood of discovering superior solutions by broadening the search during perturbation phases, albeit at a higher computational cost. A smaller $Iter_P$ amplifies the algorithm's exploratory behavior, contributing to a more thorough search of the solution space but also increasing computational demands and leading to more oscillatory convergence patterns. Conversely, selecting a larger $Amp$ facilitates wider exploration of the search space, though its effectiveness is highly contingent on the specific problem being addressed. For problems where small variations in inputs lead to significant changes in outputs, a large amplitude may not yield beneficial results, underscoring the importance of parameter tuning to align with the problem's characteristics.

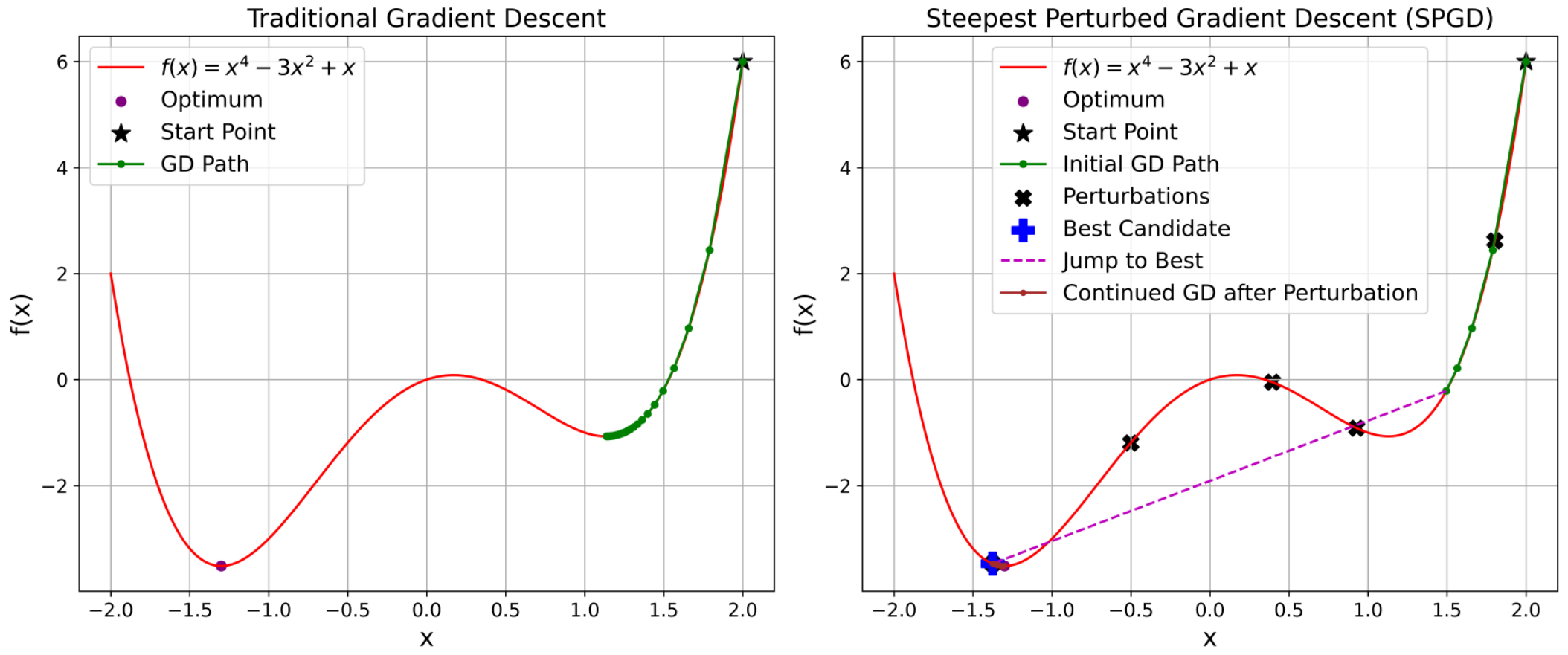

Figure 1: Optimization paths of GD and SPGD on $f(x) = x^4 - 3x^2 + x$, showing how SPGD escapes local minima and converges to the global solution.

To demonstrate the execution flow of the SPGD algorithm, we provide a visual comparison with traditional gradient descent (GD) in Figure 1, applied to a non-convex function defined by $f(x) = x^4 - 3x^2 + x$. This landscape features both a local minimum and a global minimum, offering an ideal setting to showcase the strengths of SPGD in escaping poor regions of convergence. In figure 1, GD is seen following the steepest descent path, ultimately settling at the local minimum without the ability to recover. This behavior is the characteristic of gradient-based methods in non-convex landscapes, where they are prone to getting trapped due to the absence of global information or exploration strategies.

The SPGD algorithm, on the other hand, begins similarly by following the gradient descent path, represented by the green line, for a fixed number of iterations. After $Iter_P$ iterations, the perturbation phase is triggered. At this point, SPGD generates $N_P$ candidates (depicted as black 'x') around the current position. These candidates explore the vicinity of the solution space, identifying potential points with lower function values. The most promising candidate, offering the lowest function value (shown as blue '+'), is then compared to the current position. As depicted, this candidate yields a better function value and is thus chosen, causing the algorithm to make a significant "jump" to this new point. In the figure 1, this is reflected in the trajectory of SPGD deviating from the region of slow progress of local minimum and moving toward the global one. The process then resumes with gradient descent steps leading to faster convergence to the optimal solution. As a result, SPGD not only avoids becoming stuck in local minima but also achieves convergence with fewer function evaluations compared to methods that either rely solely on gradients or purely stochastic exploration. This example highlights how SPGD's structure, consisting of gradient-based updates interleaved with perturbations, contributes to both its robustness and efficiency in solving non-convex optimization problems.

# 4 Overview of the Analysis

We present a self-contained analysis of *Steepest Perturbed Gradient Descent (SPGD)* under standard smooth nonconvex assumptions. We discuss two options:

1. **Baseline (no rollouts, our current implementation in algorithm 1):** At each perturbation round, our method samples $N_P$ candidates, scores their *immediate* objective values, selects the best immediate candidate if better than the non-perturbed (current) value, and then continues plain gradient descent (GD) until the next round. We show in this section that this variant has the same theoretical convergence rate to an approximate second-order stationary point (SOSP) as the *classical* Perturbed-GD (PGD). However, our numerical experiments described in section 5 show that this version of SPGD performs significantly better than PGD on standard optimization benchmark tests.

2. **Possible Extension to Perturbations with Rollouts:** Alternatively, one can still sample $N_P$ seeds at a perturbation round, but instead of using their immediate values, run $\tau$-step GD trajectory (a "rollout") from *each* seed, followed by the selection *after* $\tau$ steps. This strategy aligns the selection rule with the PGD escape lemma and yields a *true probability amplification* $1-(1-p_0)^{N_P}$ at the cost of additional computations.

The notation and symbols used in this section can be found in appendix B and follow those in [3].

## 4.1 Preliminaries

**Assumption 4.1** (Smoothness and lower boundedness)**.** The objective $f : \mathbb{R}^d \to \mathbb{R}$ is twice continuously differentiable and

(i) has $\ell$-Lipschitz gradient: $\|\nabla f(\mathbf{x}) - \nabla f(\mathbf{y})\|_2 \le \ell \,\|\mathbf{x}-\mathbf{y}\|_2$;

(ii) has $\rho$-Lipschitz Hessian: $\|\nabla^2 f(\mathbf{x}) - \nabla^2 f(\mathbf{y})\|_{\mathrm{op}} \le \rho \,\|\mathbf{x}-\mathbf{y}\|_2$;

(iii) is lower bounded: $f^\star := \inf_{\mathbf{x}} f(\mathbf{x}) > -\infty$.

These are standard in perturbed GD analyses; (i) yields the descent lemma - see appendix A; (ii) controls third-order effects in the saddle-escape argument [3]; (iii) lets us telescope decreases using $\Delta_f := f(\mathbf{x}_0) - f^\star$.

**Definition 4.2** (Approximate Second-Order Stationary Point (SOSP))**.** Given $\epsilon > 0$, a point $\mathbf{x}^\star$ is an $(\epsilon, \sqrt{\rho\epsilon})$-SOSP if $\|\nabla f(\mathbf{x}^\star)\|_2 \le \epsilon$ and $\lambda_{\min}(\nabla^2 f(\mathbf{x}^\star)) \ge -\sqrt{\rho\epsilon}$. The curvature threshold $\sqrt{\rho\epsilon}$ is the natural scale under standard Hessian-Lipschitz assumptions as in [3].

Throughout the analysis, points with large gradient are treated as descent region, points with small gradient but significant negative curvature as strict saddles, and the remaining small-gradient, no-negative-curvature region as the target SOSP region.

#### Shared quantities and their roles

- **Step size:** $\alpha = c/\ell$ with small universal $c \in (0, 1]$. This step size ensures $f(\mathbf{x}_{i+1}) \le f(\mathbf{x}_i) - \frac{\alpha}{2}\|\nabla f(\mathbf{x}_i)\|^2$ on non-perturbed steps - see appendix A.
- **Escape horizon:** $\tau = \Theta(\frac{\ell}{\sqrt{\rho\epsilon}} \log \frac{d}{\delta_{\text{round}}})$. With $\alpha = \Theta(1/\ell)$, $\tau$ GD steps are enough to exploit negative curvature of size $\sqrt{\rho\epsilon}$ and realize a decrease $\Omega(\epsilon^2/\ell)$ with high probability.
- **Perturbation amplitude:** $\text{Amp} = \Theta\big((\epsilon/\ell)\sqrt{\log \frac{d}{\delta_{\text{round}}}}\big)$. Radius $\epsilon/\ell$ is the scale at which a successful seed enables $\Theta(\epsilon^2/\ell)$ decrease; the $\sqrt{\log}$ factor shrinks the stuck–region volume.
- **Period between perturbations:** $\text{Iter}_P$. In the *no-rollout* case (section 4.2) $\text{Iter}_P$ is any fixed period at which a single try is attempted; In the *rollout* strategy (section 4.3) we would select $\text{Iter}_P = \tau$.
- **Number of candidates:** We select $N_P \in \mathbb{N}$ seeds per round, i.i.d. (independent and identically distributed) from the uniform distribution $\text{Unif}(\mathbb{B}(0, \text{Amp}))$.
- **Failure budgets:** total budget $\delta \in (0, 1)$; if there are $R$ rounds, we set $\delta_{\text{round}} = \delta/R$ so a union bound yields total failure $\le \delta$.

### 4.2 Baseline: Convergence *without* Rollouts (Current SPGD Implementation)

As mentioned above, our current implementation of SPGD samples $N_P$ candidates at each perturbation round, evaluates their *immediate* objectives, and selects the best immediate candidate. If the chosen candidate does not worsen the objective, we continue plain GD with the chosen candidate until the next perturbation round after $\text{Iter}_P$ steps.

Although we evaluate only the immediate objective values, this baseline mirrors the logic of the PGD perturbation step. In PGD, after injecting a single perturbation in regions where the gradient is small, one runs a $\tau$-step trajectory to test whether the perturbed point ultimately leads to a decrease in the objective beyond a fixed threshold; a perturbation is “successful” precisely when it places the iterate into a region from which GD can make progress. SPGD checks the same condition directly at the perturbation step: among the sampled candidates, we accept only if the corresponding objective value is no worse than the current iterate. However, the rounds occur after $\text{Iter}_P$ steps regardless of the gradient value.

When SPGD perturbation rounds occur near small-gradient regions, which is the same setting in which PGD perturbs, this immediate non-worsening rule serves the same purpose as PGD’s rollout test and therefore inherits the same saddle-escape guarantees. Furthermore, drawing multiple candidates increases the chance of sampling a promising direction, and never performs worse than a single PGD perturbation. We have demonstrated this on standard benchmark tests through ample experiments described in section 5.

The following theorem states the resulting guarantee for this conservative setting.

**Theorem 4.3** (**Baseline convergence without rollouts**). *Suppose Assumption 4.1 holds and suppose $\epsilon \le \ell^2/\rho$. Choose a stepsize $\alpha = c/\ell$ (small enough) with $c \in (0, 1]$, as well as* Amp *and* $\tau = \mathrm{Iter}_P$ *as described above. At each perturbation round i (separated by a period* $\mathrm{Iter}_P$*), draw $N_P$ i.i.d. perturbations $\xi^{(j)} \sim \mathrm{Unif}(\mathbb{B}(0, \mathrm{Amp}))$; set $\mathbf{x}_i^{(j)} = \mathbf{x}_i + \xi^{(j)}$ with $y_i^{(j)} = f(\mathbf{x}_i^{(j)})$; and then select the most promising seed from the $N_P$ objectives as compared with the non-perturbed objective. Next, run plain GD for $\tau = \mathrm{Iter}_P$ steps before the next round.*

*Then there exists a dimension–dependent constant $p_0 \ge 1/\mathrm{poly}(d)$ (same as in the PGD stuck–set analysis from [3]) such that the algorithm outputs an $(\epsilon, \sqrt{\rho\epsilon})$–SOSP, with probability at least $1 - \delta$, and within*

$$\tilde{O}\Big(\frac{\ell\,\Delta_f}{\epsilon^2}\Big) \text{ iterations}, \qquad \Delta_f := f(\mathbf{x}_0) - f^\star.$$

**Proof sketch:**

(i) *Descent away from stationarity:* For any non-perturbed step with $\alpha = c/\ell$, $f(\mathbf{x}_{i+1}) \le f(\mathbf{x}_i) - \frac{\alpha}{2}\|\nabla f(\mathbf{x}_i)\|^2$. Hence whenever $\|\nabla f(\mathbf{x}_i)\| \ge \epsilon$, we decrease $f$ by $\Omega(\epsilon^2/\ell)$. Summing until entering the region $\{\|\nabla f\| \le \epsilon\}$ uses at most $\tilde{O}(\ell\Delta_f/\epsilon^2)$ iterations.

(ii) *Escape from strict saddles:* If $\|\nabla f(\mathbf{x}_i)\| \le \epsilon$ and $\lambda_{\min}(\nabla^2 f(\mathbf{x}_i)) \le -\sqrt{\rho\epsilon}$, the PGD analysis of [3] shows that a random seed $\mathbf{x}_i^{(j)} = \mathbf{x}_i + \xi^{(j)}$ with radius Amp and horizon $\tau = \Theta(\ell/\sqrt{\rho\epsilon}\log(d))$ lands in an escape region with probability at least $p_0 \ge 1/\mathrm{poly}(d)$. These PGD successful seeds are guaranteed to satisfy the non-worsening condition, $f(\mathbf{x}_i^{(j)}) \le f(\mathbf{x}_i)$ (a necessary condition for escape) and are therefore accepted by SPGD. Since SPGD performs at least as well as PGD on a single sample, the per-round success probability is at least $p_0$

(iii) *Union bound and total iterations:* Let $R$ be the number of rounds, and set $\delta_{\mathrm{round}} = \delta/R$. By the union-bound (i.e., Boole's inequality [21]) across rounds we obtain a total success probability of at least $1 - \delta$. Aggregating decreases from (i)–(iii) yields the iteration bound $\tilde{O}(\ell\Delta_f/\epsilon^2)$.

**Remark 4.4** (**Multiple immediate seeds vs. single-try guarantee**). While the non-worsening rule $(f(\mathbf{x}_i^{(j)}) \le f(\mathbf{x}_i))$ is a *necessary* precursor for a successful $\tau$-step escape, the immediate objective selection alone is *not sufficient* to guarantee that the selected seed lands in the PGD escape region. Therefore, in the worst-case analysis, we rely solely on the success probability $p_0$ guaranteed by PGD's single-sample distribution. The multiple sampling factor $N_P$, which dramatically improves empirical results by acting as a strong geometric filter, is thus asymptotically nullified in the worst-case bound $\tilde{O}(\ell\Delta_f/\epsilon^2)$, leading to the same proven complexity as PGD.

**Compute cost (no rollouts).** We follow one GD trajectory between perturbations (rounds), as in PGD. Scoring $N_P$ immediate candidates at each round adds $O((T/\mathrm{Iter}_P) \cdot N_P)$ *function* evaluations; if forward passes are "cheap" compared to gradients (or the round's gradient can be reused), the overhead is modest. Crucially, since the total gradient complexity remains $O(T)$, this variant preserves the "almost dimension-free" scaling property (polylog$(d)$ dependence) of classical PGD.

### 4.3 Possible Extension: Convergence *with* Best-of-$N_P$ Rollouts

Importantly, we could make SPGD reliably find an approximate second-order stationary point in about $\tilde{O}(\ell\Delta_f/\epsilon^2)$ gradient steps by slightly changing how we use perturbation at the cost of additional computations. Specifically, instead of drawing one random perturbation and committing to it, at each "perturbation round" we could draw $N_P$ random seeds around the current point, run a short gradient descent rollout from each seed for a fixed number of steps $\tau$, and then continue the algorithm from whichever seed achieved the lowest objective value after those $\tau$ steps. Between perturbation rounds we would just carry out normal gradient descent so we can only spend about $\tilde{O}(\ell\Delta_f/\epsilon^2)$ iterations before the gradient becomes small.

In the neighborhood of a region where the gradient is small but there is a strongly negative curvature direction (a strict saddle), prior analyses [3] show that a single random perturbation followed by $\tau$ gradient steps has a nontrivial probability (at least $1/\text{poly}(d)$) of making a significant decrease in the function value. By running $N_P$ independent rollouts in parallel and picking the best one, we can amplify this success probability to be very close to one in each round. With appropriate choices of step size, perturbation radius, rollout length $\tau$, and number of seeds $N_P$ tied to the smoothness and curvature parameters of the function, we could guarantee that each time we approach a strict saddle we escape it quickly with high probability.

Because the theoretical escape event is defined in terms of improvement after $\tau$ steps, selecting the best candidate at time $i + \tau$ ensures that whenever a seed produces the escape, the decrease corresponding to the chosen seed is at least as large as the guaranteed decrease from that event.

### 4.3.1 Time and Compute Complexity for the Best-of-$N_P$ Rollouts Case

Denote by $T = \tilde{O}(\ell\Delta_f/\epsilon^2)$ the total number of GD iterations until termination. Let $R \asymp T/\text{Iter}_P$ be the number of perturbation rounds, and recall $\tau = \text{Iter}_P$.

**Gradient evaluations** In this scenario, we would always execute $T$ gradient steps along the selected path. In addition, at each round we roll out $N_P - 1$ extra seeds for $\tau$ steps that are discarded after selection. Thus the total number of gradient evaluations would be

$$T \;+\; R\,(N_P - 1)\,\tau \;=\; T \;+\; \frac{T}{\text{Iter}_P}\,(N_P - 1)\,\text{Iter}_P \;=\; O(N_P\,T).$$

Hence amplification would trade a linear (in $N_P$) increase in computational cost for a much smaller per-round failure probability.

**Function evaluations** If the implementation can reuse the round's gradient $\nabla f(\mathbf{x}_i)$ for all seeds and only needs forward evaluations to score candidates after rollout, the additional forward cost per round is $O(N_P\,\tau)$ (often cheaper than $N_P\,\tau$ full gradient steps). Otherwise, the cost scales comparably to the gradient bound above.

**Asymptotic Iteration Order For Both Strategies** These $N_P$ rollouts would buy us a higher per-round reliability (smaller failure probability, so better constants/polylogs[3], but they would not increase the per-success decrease nor the per-step decrease in large-gradient regions. That is why both approaches are expected to have the same asymptotic iteration order $\tilde{O}(\ell\,\Delta_f/\epsilon^2)$.

## 4.4 Comparison with PGD

Both variants match PGD's iteration complexity $\tilde{O}(\ell\Delta_f/\epsilon^2)$ to reach an $(\epsilon, \sqrt{\rho\epsilon})$-SOSP. We note that our extensive experiments described in section 5 demonstrate that the no-rollout SPGD of section 4.2 performs significantly better than PGD on standard benchmark functions due to a better selection of candidates from the neighborhood. A careful analysis between the tradeoff between efficiency and accuracy of the rollout mechanism as well as its implementation are outside the scope of this paper, but are certainly a logical next step. This is so because the rollout variant of section 4.3 trades, as already mentioned, an $O(N_P)$ compute factor for an exponentially reduced per-round failure.

[3]We use the notation $\text{polylog}(n)$ to denote a polynomial in $\log n$, i.e., $\text{polylog}(n) = (\log n)^{O(1)}$. Accordingly, we write $\widetilde{O}(f(n))$ for $\widetilde{O}(f(n)) := O\big(f(n)\ \text{polylog}(n)\big)$; see, e.g., [22].

# 5 Numerical Results

We present here a thorough evaluation of the proposed Steepest Perturbed Gradient Descent (SPGD) algorithm, comparing its performance against several established optimization methods. The comparison includes traditional gradient descent (GD), Perturbed Gradient Descent (PGD), MATLAB $fmincon$ function, which is a versatile solver for constrained optimization problems [23], and the $fminunc$ function, which is tailored to unconstrained optimization problems [24], along with Simulated Annealing (SA) [25], and Bayesian Optimization (BO) [26].

Our initial analysis is conducted through the lens of four challenging 2D benchmark functions, selected for their known difficulties and relevance in assessing optimization algorithms' efficacy. These test functions are recognized benchmarks within the optimization community, providing a diverse set of landscapes to evaluate each algorithm's ability to navigate complex, non-convex, and potentially deceptive optimization spaces [27]. For each test function, we apply $fmincon$, Simulated Annealing, traditional gradient descent, PGD, and SPGD, meticulously recording and analyzing the results. The SPGD and PGD algorithms were each fine-tuned independently to ensure optimal performance while maintaining a fair basis for comparison. Due to the high computational cost of Bayesian Optimization (BO), the maximum number of function evaluations for BO is capped at 100 to ensure reasonable execution time across benchmark functions. For the $fmincon$ function, we use MATLAB's default *interior-point* algorithm, while the $fminunc$ function is configured to use the *trust-region* method, which is well-suited for smooth, unconstrained problems. In both cases, the gradient of the objective function is explicitly provided to guide the optimization process more efficiently.

Key performance indicators include the accuracy of the solution, measured by the proximity to the known global optimum [28]; the computational efficiency, quantified by the number of function evaluations and CPU execution time. For each test function, both a 3D and top-view surface plot of optimization trajectory visualization are provided to aid in understanding each algorithm's optimization landscape and behavior. These visualizations illustrate how optimization paths evolve over complex response surfaces and help highlight differences in convergence dynamics. Simulated Annealing (SA) and Bayesian Optimization (BO) are excluded from these visualizations, as their probabilistic sampling strategies tend to densely populate the landscape, obscuring the trajectories of other algorithms and reducing the overall interpretability of the plots. The source code for the SPGD algorithm, along with comparative analyses against methods discussed in this paper using additional 2D challenging test functions, are publicly accessible on GitHub[4].

### Test function 1

The MATLAB Peaks function [29] presents a formidable challenge for optimization algorithms due to its intricate landscape, which features one global minimum, multiple local minima, a saddle point, and extensive flat regions. This complexity makes the Peaks function a critical benchmark for assessing the capabilities of optimization techniques, particularly those based on gradient descent. Traditional gradient descent methods often struggle with such landscapes, as they can easily become trapped in local minima or stall in flat areas, failing to make significant progress towards the global optimum [30]. The mathematical expression defining the Peak test function is given as follows:

$$f(x, y) = 3\,\mathrm{e}^{-(y+1)^2-x^2}\,(x-1)^2 - \frac{\mathrm{e}^{-(x+1)^2-y^2}}{3} + \mathrm{e}^{-x^2-y^2}\left(10\,x^3 - 2\,x + 10\,y^5\right) \quad (2)$$

It has a global minimum point located at $x = 0.2283$, $y = -1.6256$ with an optimal function value of $f(x^*) = -6.5511$. The initial condition is chosen randomly to be $(-2.81, -1.47)$, and the $Amp$ is set to 2.5. Figure 2a and 2b illustrate the 3D view and top view of the optimization trajectory across the Peaks function surface. The total number of function evaluations, the converged optimal value, and CPU

[4]Source code and comparisons available at: `https://github.com/Amir-M-Vahedi/SPGD-Benchmark-Functions`

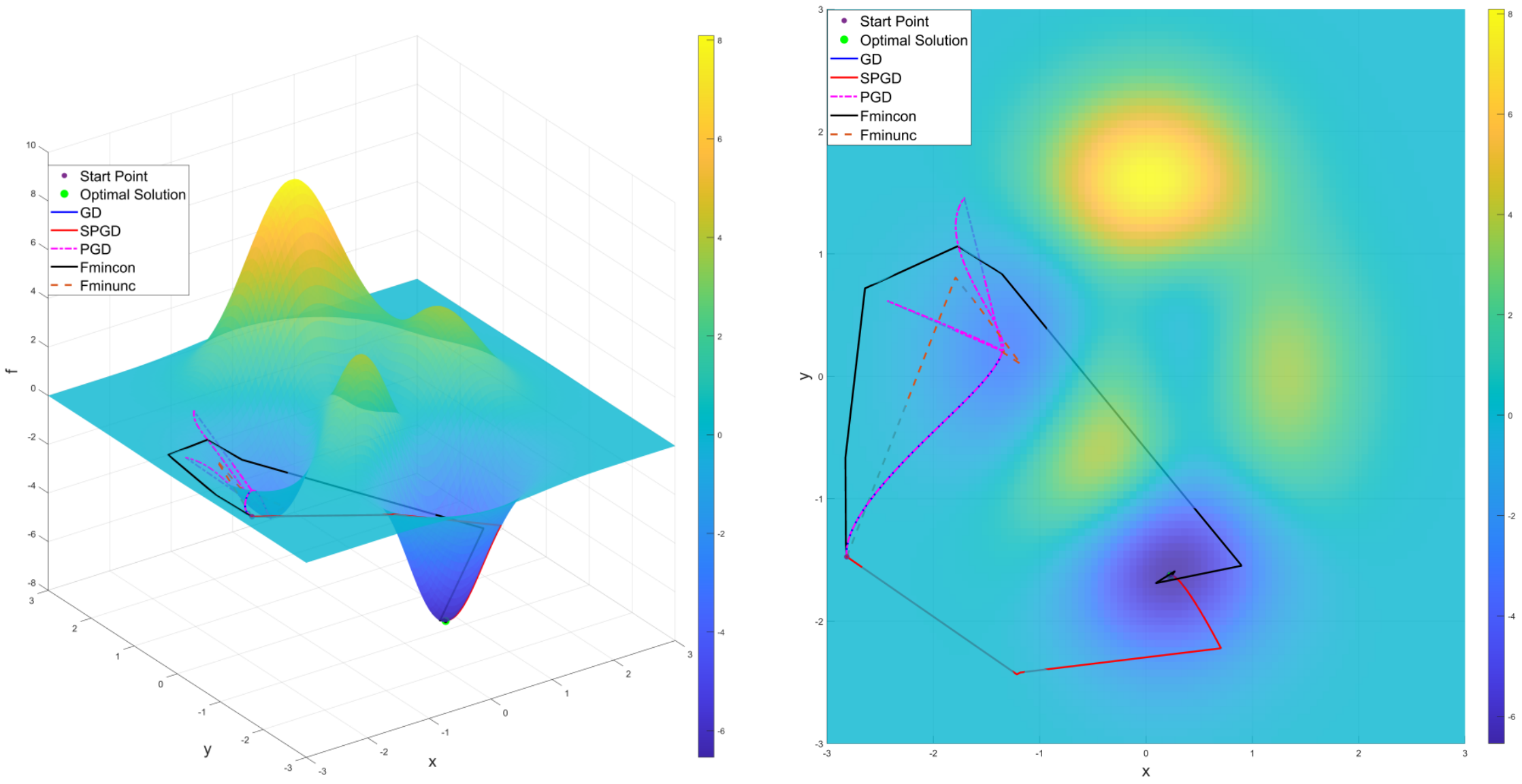

(a) 3D view of the optimization trajectories.

(b) Top view of the optimization trajectories.

Figure 2: Visualization of optimization trajectories for the Peaks function.

execution time for different methods are given in Table 1. Based on the results depicted in Figures 2, and performance metrics in Table 1, it is evident that the GD, PGD, and $fminunc$ algorithms become trapped in local minima. In contrast, the $fmincon$, SA, BO, and SPGD algorithms successfully converge to the global optimum. Among these three, SPGD demonstrates the lowest computational cost. Notably, despite the $fmincon$ and BO method having fewer function evaluations, their CPU times are more than 25 and 2000 times greater than that of the SPGD algorithm.

Table 1: Peaks function Performance

| **Algorithm** | **Total Fun. Evaluations** | $f(x^*) = -6.5511$ | **CPU Time[ms]** |
|---|---|---|---|
| GD | 1472 | -3.0498 | *3.12 |
| PGD | 1599 | -3.0498 | *3.88 |
| $fminunc$ | 10 | -3.0498 | *23.25 |
| $fmincon$ | 60 | **-6.5511** | 57.37 |
| SA | 1341 | **-6.5511** | 117.8327 |
| BO | 100 | **-6.5510** | 4415.7 |
| **SPGD** | 274 | **-6.5511** | 2.03 |

### Test function 2

The Ackley function is a well-known non-convex optimization benchmark that poses a significant challenge to optimization algorithms, particularly due to its deceptive landscape characterized by a global optimum surrounded by a multitude of local minima [31]. This function is specifically designed to test the ability of optimization methods to escape local minima and efficiently search for the global optimum in a complex, multidimensional space. The Ackley function's landscape features a large number of local minima leading towards the global minimum, making it an exemplary test case for evaluating the robustness and effectiveness of algorithms against the risk of premature convergence. The global minimum of the Ackley function is located at the origin ($x = 0$, $y = 0$), with an optimal function value of zero ($f(x^*) = 0$), which further serves as a clear target for optimization efforts. The formula representing the Ackley test function

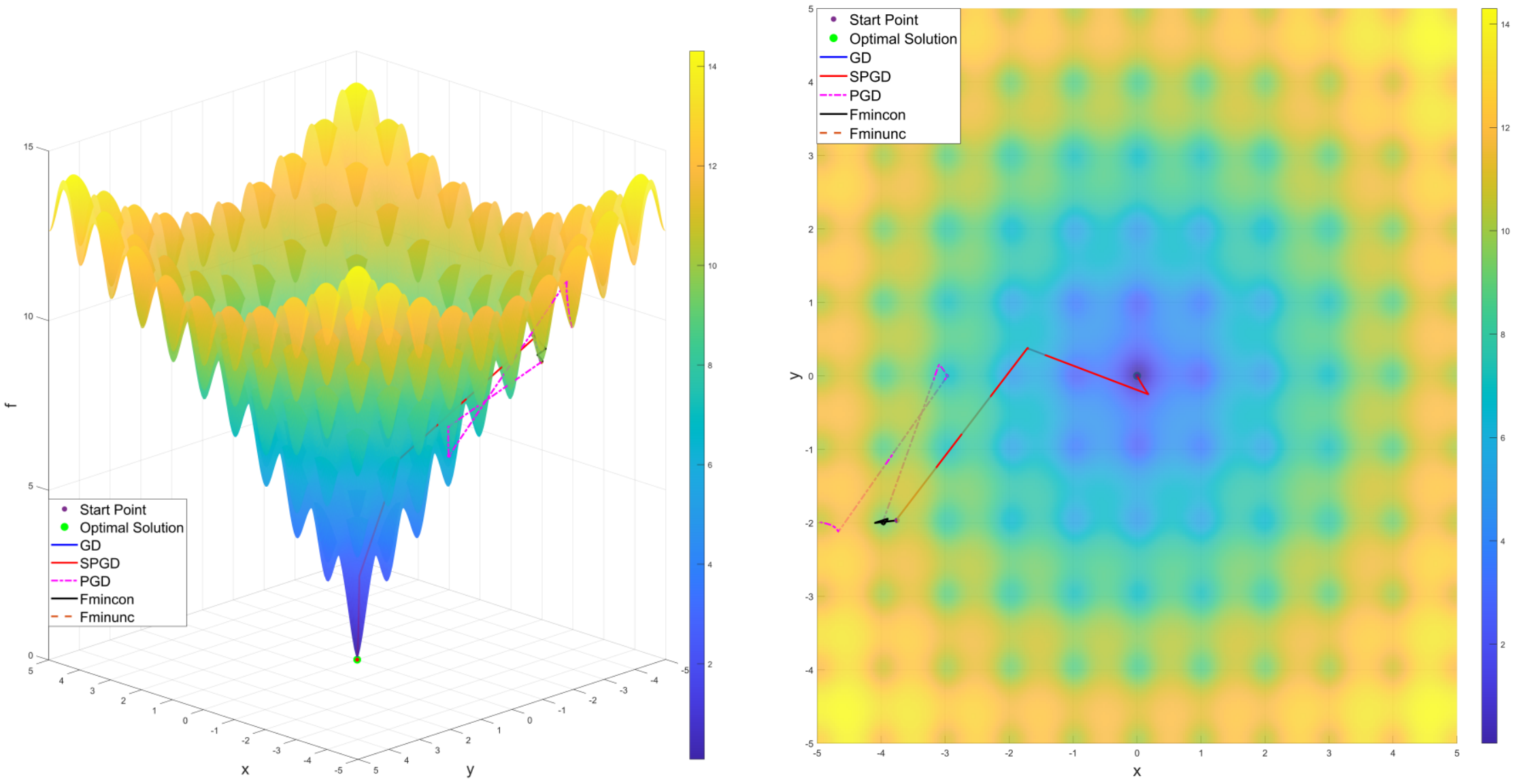

(a) 3D view of the optimization trajectories. (b) Top view of the optimization trajectories.

Figure 3: Visualization of optimization trajectories for the Ackley function.

is articulated below:

$$f(x,y) = -20\exp\left(-0.2\sqrt{\frac{1}{2}(x^2+y^2)}\right) - \exp\left(\frac{1}{2}(\cos(2\pi x)+\cos(2\pi y))\right) + 20 + e \tag{3}$$

The initial condition is chosen randomly to be $(-3.75, -1.96)$, and the $Amp$ is set to 2.5. Figure 3a and 3b illustrate the 3D view and top view of optimization trajectory across the Ackley function surface. The performance comparisons are given in Table 2. Taking into account the data presented in the mentioned figures and table, analysis of the Ackley test function reveals that the GD, PGD, $fminunc$, and $fmincon$ methods became ensnared in local minima. In contrast, only the SA, BO, and SPGD algorithms successfully navigated to the global solution. Among these, SPGD not only achieved convergence with greater precision, approaching closer to the global optimum, but also demonstrated a computational speed, with the CPU execution time being about 13 (SA) and 359 (BO) times faster than its counterparts.

Table 2: Ackley function Performance

| **Algorithm** | **Total Fun. Evaluations** | $f(x^*) = 0$ | **CPU Time[ms]** |
|---|---|---|---|
| GD | 327 | 9.3530 | *2.01 |
| PGD | 477 | 6.8826 | *2.02 |
| $fminunc$ | 8 | 9.3530 | *32.20 |
| $fmincon$ | 24 | 9.3530 | *75.13 |
| SA | 504 | **2.13e-4** | 40.63 |
| BO | 100 | **0.0213** | 4670.0 |
| **SPGD** | 1501 | **4.81e-4** | 3.62 |

### Test function 3

The Easom function stands as a notable unimodal steep ridge [27] test function within the realm of optimization, particularly distinguished by its singular global optimum that resides in an extensive flat area. This flat region is characterized by minimal gradient variations, presenting a unique challenge for optimization algorithms, especially those reliant on gradient information to navigate the search space. The function is defined over a domain of $(-100, 100)$ for both $x$ and $y$ dimensions, emphasizing the necessity for optimization techniques to efficiently explore large search areas to locate the optimum [32]. The significance of the Easom function as a test scenario with simple mathematical formulation lies in its ability to simulate real-world optimization problems where the solution space is largely homogeneous, yet contains a singular, critical point of interest. This function tests the exploration strategies of algorithms, challenging them to avoid the pitfalls of vast non-informative regions. It emphasizes the importance of balance between exploration and exploitation, as effective optimization methods must not only navigate vast spaces efficiently but also recognize and converge to the global optimum with high precision. Mathematically, the Easom function's global optimum is uniquely situated at $(x = \pi,\ y = \pi)$, where it attains a value of $f(x^*) = -1$. The formula of the Easom test function is provided below:

$$f(x, y) = -\cos(x)\cos(y)\exp\Big(-\Big((x-\pi)^2 + (y-\pi)^2\Big)\Big) \tag{4}$$

The initial condition is chosen randomly to be $(69.33, 12.23)$, and the $Amp$ is set to 5. Figure 4a and 4b illustrate the 3D view and top view of the optimization trajectory across the Easom function surface. The performance comparisons are given in Table 3. Reflecting on the performance metrics for the Easom test function, it is evident that only the SPGD algorithm successfully pinpointed the global optimum. As anticipated, GD was hindered in its progression by the minimal gradient values inherent to the function's extensive flat regions. Similarly, both GD and the $fmincon$ method failed to escape these flat expanses, effectively becoming ensnared within them. Among the competing methods, only SA and BO managed to navigate towards a more favorable outcome, yet they fell short of achieving convergence to the global optimum, underscoring the distinctive effectiveness of SPGD in this scenario.

Table 3: Easom function Performance

| **Algorithm** | **Total Fun. Evaluations** | $f(x^*) = -1$ | **CPU Time[ms]** |
|---|---|---|---|
| GD | 1 | 0 | *0.08 |
| PGD | 2021 | 0 | *3.28 |
| $fminunc$ | 1 | 0 | *58.35 |
| $fmincon$ | 1 | 0 | *83.84 |
| SA | 1009 | -3.38$e$-160 | *69.60 |
| BO | 100 | -6.22$e$-215 | *6866.9 |
| **SPGD** | 6001 | **-1** | 7.45 |

**Test function 4**

The Levy Function No. 13, characterized by its multimodality and non-convexity, presents a unique challenge for optimization algorithms with its single global optimum amidst a noisy, periodic distribution of local minima. This function tests an algorithm's precision in distinguishing the global optimum from numerous suboptimal states, a key trait for solving complex real-world problems. It serves as a critical benchmark for evaluating the balance between exploration and exploitation in optimization techniques, underscoring its significance in both theoretical and practical applications. The global optimum of this function is strategically located at $(x = 1,\ y = 1)$, where it attains a value of $f(x^*) = 0$. The expression for the Levy Function No. 13 is detailed below [33]:

$$f(x, y) = \sin^2(3\pi x) + (x-1)^2\Big(1 + \sin^2(3\pi y)\Big) + (y-1)^2\Big(1 + \sin^2(2\pi y)\Big) \tag{5}$$

The initial condition is chosen randomly to be $(-3.75, -1.96)$, and the $Amp$ is set to 2.5. Figure 5a and 5b illustrate the 3D view and top view of the optimization trajectory across the Levy Function No. 13

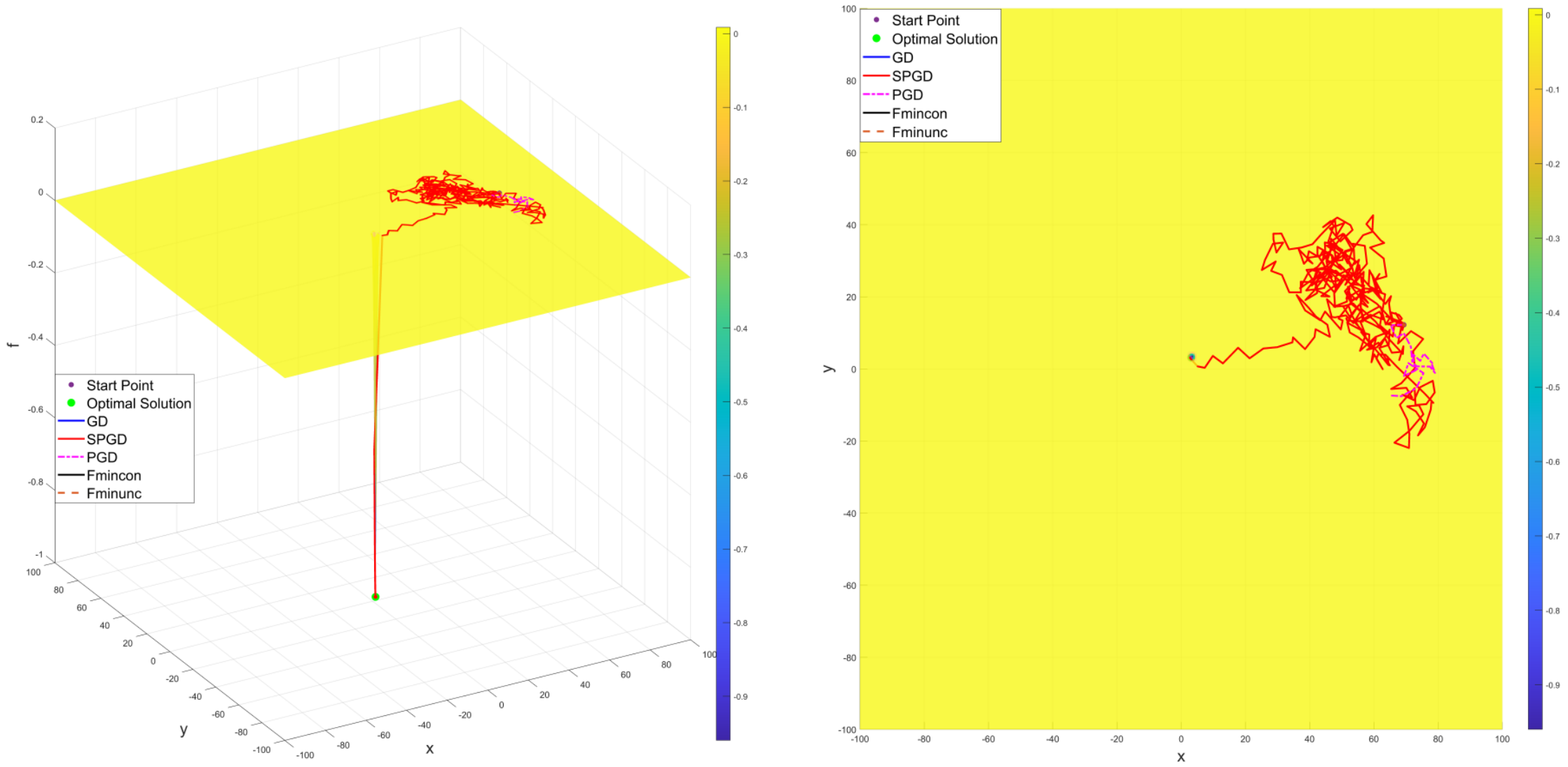

(a) 3D view of the optimization trajectories. (b) Top view of the optimization trajectories.

Figure 4: Visualization of optimization trajectories for the Easom function.

surface. The performance comparisons are given in Table 4. Based on the performance analysis for this test function, the GD, PGD, *fminunc*, and *fmincon* methods were unable to find the global optimum, getting stuck in local minima instead. Notably, *fmincon* settled in a particularly poor local minimum. In contrast, SA, BO, and SPGD successfully navigated to the global optimum. However, SPGD distinguished itself by achieving a more accurate solution, requiring fewer function evaluations than SA, and demonstrating faster CPU execution time compared to SA and BO.

Table 4: Levy function N. 13 Performance

| **Algorithm** | **Total Fun. Evaluations** | $f(x^*) = 0$ | **CPU Time[ms]** |
|---|---|---|---|
| GD | 2001 | 6.2915 | *3.58 |
| PGD | 2001 | 6.2915 | *2.58 |
| *fminunc* | 9 | 14.3717 | *20.59 |
| *fmincon* | 20 | 30.5009 | *53.32 |
| SA | 2018 | **6.78*e*-7** | 89.61 |
| BO | 100 | **0.0086** | 5241.7 |
| **SPGD** | 1760 | **2.45*e*-13** | 5.02 |

## Robustness Evaluation

To evaluate the robustness and statistical reliability of the SPGD algorithm, each benchmark function was tested using 30 independent trials with randomly sampled initial points. All optimization algorithms were provided with the same lower and upper bounds defining the feasible search space, and each method was fine-tuned independently to ensure its best individual performance. The following performance criteria are reported in Tables 5–8:

To evaluate the robustness of the SPGD algorithm, each test function was subjected to 30 independent trials using randomly sampled starting points. All optimization algorithms were provided with the same lower and upper bounds defining the feasible search space. Each algorithm was fine-tuned independently to ensure their optimal individual performance, and consistent parameter settings were applied across all

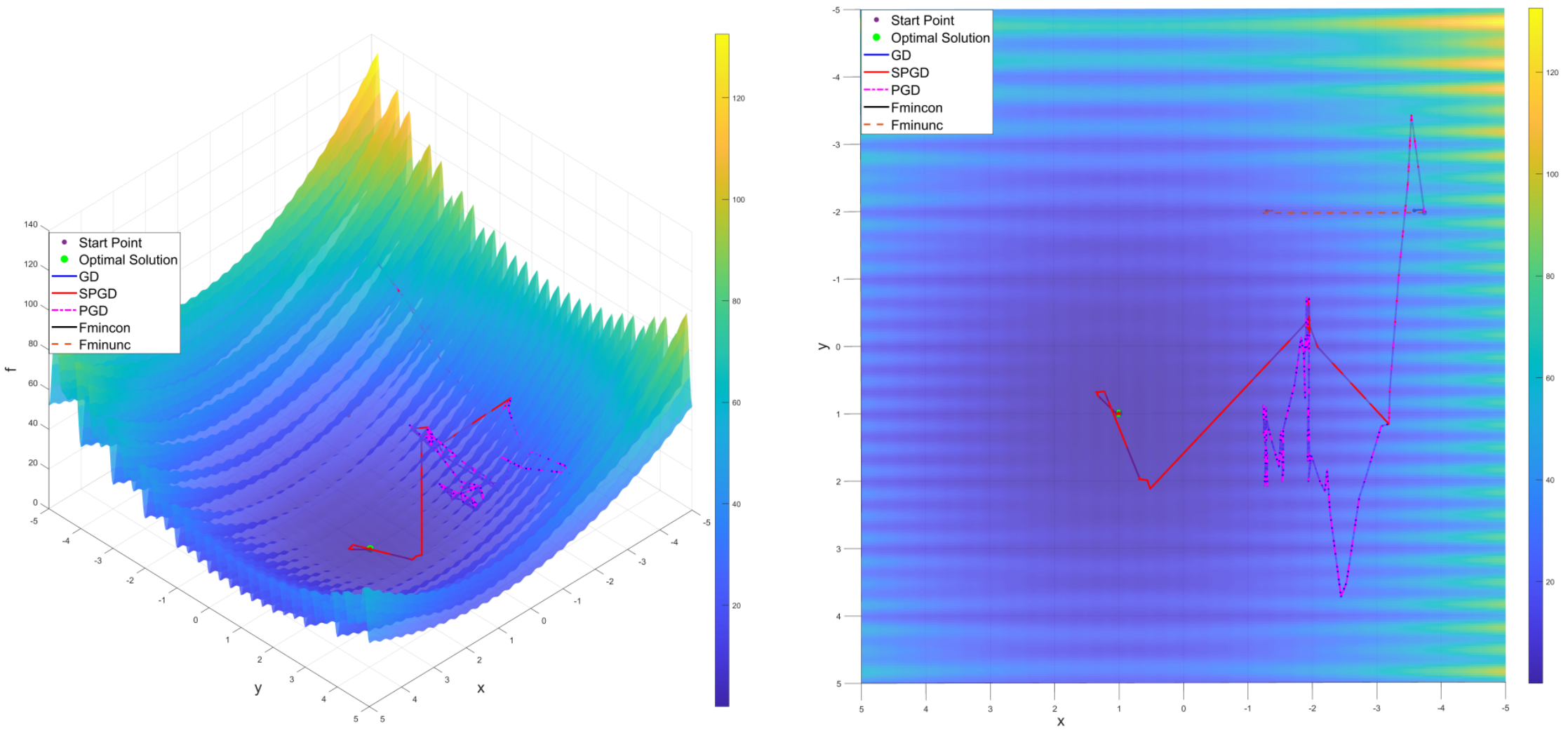

(a) 3D view of the optimization trajectories. (b) Top view of the optimization trajectories.

Figure 5: Visualization of optimization trajectories for the Levy function N. 13.

trials to maintain fairness. The following performance criteria are reported in Tables 5–8:

- **ConvergedRuns**. This metric reports how many out of the 30 randomized trials a given method successfully converged to the global optimum within a tolerance of $10^{-6}$.
- **Fval Improvement (%)** and **Time Improvement (%)**. These metrics quantify the relative difference between each method and SPGD, expressed as a percentage. For a given method $i$, the percentage difference in objective value is computed as

$$\text{Improvement}_{\text{fval}}(i) = 100 \times \frac{\bar{f}(i) - \bar{f}_{\text{SPGD}}}{|\bar{f}(i)|}, \tag{6}$$

  where $\bar{f}(i)$ is the mean objective value of method $i$ over 30 trials, and $\bar{f}_{\text{SPGD}}$ is the corresponding mean for SPGD. A positive value indicates that method $i$ performs worse than SPGD (higher objective value), while a negative value indicates better performance. The same expression is used for CPU time by replacing objective values with average execution times.

- **Closer (%)**. This metric quantifies how much closer method $i$ is to the global optimum compared to SPGD. Let $f^*$ denote the known global minimum of the benchmark function. Define the mean distance to the optimum for method $i$ as

$$d(i) = |\bar{f}(i) - f^*|. \tag{7}$$

  The closeness metric is then computed as

$$\text{Closer}(i) = 100 \times \frac{d(i) - d_{\text{SPGD}}}{d(i)}, \tag{8}$$

  where $d_{\text{SPGD}}$ is the distance for SPGD. Positive values indicate that method $i$ is farther from the optimum than SPGD, while negative values indicate that it is closer. When the comparing method also attains the exact global optimum (i.e., $d(i) = 0$), this metric is not defined and is reported as “N/A”.

For the Peaks function, SPGD successfully converged to the global optimum in all 30 trials. Only Bayesian Optimization (BO) matched this convergence count, with Simulated Annealing (SA) achieving 29 out of 30. However, SPGD accomplished this with significantly lower computational cost, as reflected in the *Time Improvement%* column of Table 5, demonstrating superior efficiency.

In the case of the Ackley function, SPGD was the only algorithm to consistently converge to the global solution across all runs. Although GD and PGD had lower average execution times, they only succeeded in 3 and 5 out of 30 runs respectively, making them less competitive. SPGD outperformed all remaining methods in terms of average speed and reliability over 30 trials.

For the Easom function, none of the baseline algorithms found the global optimum in any run. SPGD was the only method to successfully reach the global solution in all 30 trials, highlighting its robustness in highly deceptive landscapes.

In the Levi function N.13, SPGD again demonstrated the highest reliability with 30 successful runs out of 30. BO and SA achieved 24 successful runs each, while other algorithms failed to find the global optimum in any trial. SPGD also demonstrated strong efficiency, with an average *Time Improvement%* of 99.91 over BO and 90.02 over SA.

Table 5: Average performance comparison for Peaks function over 30 random starting points

| **Algorithm** | **ConvergedRuns** | **Fval Improvement %** | **Time Improvement %** | **Closer %** |
|---|---|---|---|---|
| GD | 9 | 135.63 | -5.43 | 100.00 |
| PGD | 9 | 127.32 | -10.01 | 100.00 |
| BayesOpt | 30 | 0.00 | 99.99 | N/A |
| SA | 29 | 0.04 | 98.90 | 98.71 |
| Fminunc | 7 | 206.33 | 81.57 | 100.00 |
| Fmincon | 8 | 137.05 | 93.10 | 100.00 |
| SPGD | 30 | | | |

Table 6: Average performance comparison for Ackley function over 30 random starting points

| **Algorithm** | **ConvergedRuns** | **Fval Improvement %** | **Time Improvement %** | **Closer %** |
|---|---|---|---|---|
| GD | 3 | 99.97 | -1240.80 | 99.97 |
| PGD | 5 | 99.97 | -1079.38 | 99.97 |
| BayesOpt | 17 | 83.04 | 99.96 | 83.04 |
| SA | 27 | 99.15 | 85.82 | 99.15 |
| Fminunc | 5 | 99.97 | 50.29 | 99.97 |
| Fmincon | 8 | 99.96 | 83.92 | 99.96 |
| SPGD | 30 | | | |

The challenges presented by these test functions, including their rugged landscapes and deceptive local minima, contain features that bear resemblance to those encountered in the energy landscape of protein folding. This complex biological process is characterized by a similarly intricate energy landscape that features multiple local optima (kinetic traps), rugged terrain, and steep energy barriers (sharp valleys and hills) [34–37].

The SPGD algorithm's performance on these test functions suggests its potential utility in addressing the complex optimization problems inherent in protein folding. By adeptly navigating through challenging landscapes to find global or near-global optima, SPGD could significantly contribute to bioinformatics and molecular biology by optimizing protein structures to understand their function and interactions more accurately.

Table 7: Average performance comparison for Easom function over 30 random starting points

| Algorithm | ConvergedRuns | Fval Improvement % | Time Improvement % | Closer % |
|---|---|---|---|---|
| GD | 0 | N/A | -79148.81 | 100.00 |
| PGD | 0 | 444556889.89 | -1822.73 | 100.00 |
| BayesOpt | 0 | 4137.02 | 99.97 | 100.00 |
| SA | 0 | 8841.39 | 78.76 | 100.00 |
| Fminunc | 0 | N/A | -724.67 | 100.00 |
| Fmincon | 0 | N/A | -277.97 | 100.00 |
| SPGD | 30 | | | |

Table 8: Average performance comparison for Levi function N.13 over 30 random starting points

| Algorithm | ConvergedRuns | Fval Improvement % | Time Improvement % | Closer % |
|---|---|---|---|---|
| GD | 0 | 100.00 | -1761.37 | 100.00 |
| PGD | 0 | 100.00 | -1462.19 | 100.00 |
| BayesOpt | 24 | 100.00 | 99.91 | 100.00 |
| SA | 24 | 100.00 | 90.02 | 100.00 |
| Fminunc | 0 | 100.00 | -140.21 | 100.00 |
| Fmincon | 0 | 100.00 | 34.57 | 100.00 |
| SPGD | 30 | | | |

This analogy not only highlights the broader applicability of SPGD but also underscores the importance of developing robust optimization techniques that can effectively deal with the complexities of both mathematical functions and biological systems [38, 39].

Expanding our investigation beyond conventional 2D test functions, we also apply our algorithm (SPGD) to a 3D component packing problem, a task distinguished by its NP-hard classification [40, 41]. This problem introduces a unique set of challenges, including flat area saddle points and local optima, that further test the robustness of our approach against traditional gradient descent and simulated annealing methods.

### 3D Component Packing Problem

In the 3D component packing problem, we focus on arranging 3D objects with arbitrary shapes as compactly as possible without collision, akin to a simplified version of the interconnected systems with physical interactions (SPI2) problem but without considering the routing interconnections between objects [42]. Optimization methods often face challenges in this landscape, such as getting trapped in local minima or stalling in flat areas, thus failing to advance significantly towards the global optimum. The non-convex nature of the objective function, characterized by multiple local optima and saddle points, poses substantial challenges to any standard optimization technique. Nonetheless, our SPGD algorithm, which integrates randomized perturbations, is tailored to navigate these complex landscapes more effectively, demonstrating its adaptability and enhanced performance compared to conventional techniques.

Our 3D packing scenarios presented here are specialized instances of those tackled by SPI2-F – a novel and more general packing and layout optimization presented in [43], that performs both packing and layout optimization of complex interconnected systems in a multi-physics environment. The specialized scenarios presented below have been chosen to include cases that have known global optima. We note that an efficient packing method based on Fast Fourier Transforms[5] was introduced recently in [48], which restricts

[5] Fast Fourier Transforms have been previously shown to offer an elegant and efficient approach to compute collisions and

the orientation of the objects to an axis-alignment and hence allows rotations in 90-degree increments. By contrast, our method allows arbitrary rotations and alignments in space.

In the 3D component packing problem, our primary objectives are twofold: minimize the volume of the bounding box containing the components ($V_b$) while avoiding collisions between the components [43]. Therefore, we define the mathematical expression of the objective function as follows:

$$f = w_b \times V_b - w_c \times \log(\epsilon + \min(\text{dist})) \tag{9}$$

where $w_b = 20$ represents the weight associated with the bounding box volume, $w_c$ =1$e$-4 is the weight related to collision avoidance, $\epsilon$ =1$e$-5 is a small value to avoid singularity, and min(dist) denotes the minimum distance between the spheres of different components.

The complexity and high dimensionality of this problem are underscored by the fact that each object in our example consists of $num_{sphere} = 100$ spheres, and each component is controlled by six variables – three for displacement and three for orientation. The problem also incorporates constraints related to collision avoidance. To effectively navigate the highly non-convex and constrained space of the component packing problem, our approach involves tailored adaptations to the perturbation mechanism used in the Steepest Perturbed Gradient Descent (SPGD) algorithm. Perturbations are applied separately to the components' displacement and orientation, ensuring a thorough optimization of both aspects of component placement.

In the early iterations, we enhance the exploration and facilitate the escape from suboptimal solutions by accepting solutions with worse volume outcomes by a prescribed factor. This acceptance factor decreases in a linear profile over the iterations until it reaches 1.0, at which point the algorithm only accepts new solutions that have the same or lower volume, thus refining the search towards the most compact configurations. Additionally, the amplitude of the perturbations for both displacement and orientation is controlled through a lower-bounded linear profile, which ensures that perturbations decrease in magnitude as the optimization process progresses, aligning more closely with the finer adjustments needed as the solution space is narrowed down. To further optimize the perturbation process and avoid ineffective perturbations, especially in cases where objects are too close to each other to allow for meaningful spatial adjustments, the frequency of perturbations is reduced using a lower bounded linear profile. This adaptive frequency adjustment helps prevent unnecessary computational expenditure on perturbations that are unlikely to be accepted due to collision constraints. These strategic adaptations enable SPGD to more effectively handle the complexities of packing diverse objects into a constrained space, making it robust against the challenges posed by the non-convex nature of the problem.

The implementation of this algorithm is carried out in Python using the PyTorch framework, which leverages CUDA for accelerated computation on GPUs. This setup allows for substantial improvements in computational efficiency, essential for managing the high-dimensional space of this problem. Using `torch.autograd`, we automatically compute the partial derivatives of the loss function with respect to the displacement and orientation vectors of each component. This gradient information is then used to update the component positions and orientations according to the update rule of gradient descent (1), akin to methods typically employed in deep learning optimizations. To further enhance the exploration capabilities of the optimization process, the sequence of component perturbations is shuffled in each iteration, promoting a more robust search through the solution space. In evaluating the effectiveness of the SPGD algorithm, we conducted a comparative analysis with the traditional Gradient Descent (GD) method across a series of increasingly complex packing scenarios.

The scenarios were designed to assess both algorithms under various conditions, ranging from uniform object sizes to irregular and diverse shapes, thereby testing their adaptability and efficiency in real-world packing challenges. Our implementation of Simulated Annealing diverged rather than converged, particularly in complex scenarios. This divergence can largely be attributed to the restrictive collision constraints integrated within the objective function (9), which prevent objects from moving through each other. Unlike the approach taken in reference [41], where collision constraints were relaxed and followed penetrations as well as shape complementarity [44–47].

by refinement steps, our implementation maintained these constraints, leading to no evident signs of convergence and indicating the unsuitability of Simulated Annealing for these applications.

Moreover, Perturbed Gradient Descent (PGD) was not utilized in the 3D packing problem comparison. The reason for this is twofold: firstly, the norm of the partial derivative vector in this problem setting does not approach zero due to the direct inclusion of collision constraints within the objective function. Secondly, the primary cause for algorithm termination is often the occurrence of collisions between objects, which deviates from the typical operational premise of PGD. Additionally, PGD's poor performance in separate 2D benchmark functions, which feature complex and challenging loss landscapes, further illustrates its limitations in navigating complicated optimization scenarios. This combination of factors reaffirms the decision to exclude PGD from the comparative analysis in our 3D packing problem.

## 5.1 Initial Configuration and Setup

Before delving into the comparative results, it is essential to note that both the SPGD and GD algorithms were initiated from the same configuration in each scenario to ensure a fair comparison. The initial setup involved distributing the objects well within the 3D space, providing sufficient free space around each object to avoid immediate collisions. Furthermore, the orientations of the objects were randomly chosen, introducing additional complexity and ensuring that the problem remained challenging for the optimizers. This initialization strategy was crucial for testing the algorithms' abilities to effectively explore and optimize from a non-advantageous starting point.

## 5.2 Experimental Scenarios and Results

The following scenarios were considered for the comparison:

- **Scenario 1: Four identical rectangular boxes.** In this case, the global optimum is analytically known. Since all objects are homogeneous rectangular boxes, the minimum-volume packing configuration corresponds to placing the rectangular boxes without gaps and with identical orientation, forming an axis-aligned rectangular block. For four rectangular boxes, this arrangement results in either a $1 \times 4$ or a $2 \times 2$ layout, both achieving the same optimal bounding-box volume.

- **Scenario 2: Eight identical rectangular boxes.** Similar to Scenario 1, the global optimum is known analytically. Any axis-aligned arrangement that packs all eight rectangular boxes into a gap-free block attains the minimum possible volume. Examples include $1 \times 8$, $2 \times 4$, or $2 \times 2 \times 2$ layouts, all of which represent globally optimal solutions.

- **Scenario 3: Eight rectangular boxes of varying sizes.** In this non-uniform configuration, no closed-form global optimum is known. The objective is to evaluate each algorithm's ability to navigate a heterogeneous packing space where optimality cannot be verified analytically.

- **Scenario 4: Eight complex-shaped objects.** This scenario includes industrial parts such as gears, hooks, rivets. Due to the geometric complexity and lack of symmetry, no analytical global optimum exists. This setup tests the heuristic and exploratory capabilities of the algorithms, representing an industrial challenge with an unknown optimal packing configuration.

## 5.3 Analysis of Scenario 1: Four Identical rectangular boxes

In Scenario 1, the initial configuration of the four identical rectangular boxes is depicted in Figure 6. This setup was designed to test each algorithm's ability to navigate a relatively simple scenario where the global optimum involves aligning all rectangular boxes in a compact configuration. The results of the final configurations found by the GD and SPGD algorithms are illustrated in Figure 7, showing both Gradient Descent and Steepest Perturbed Gradient Descent results side by side.

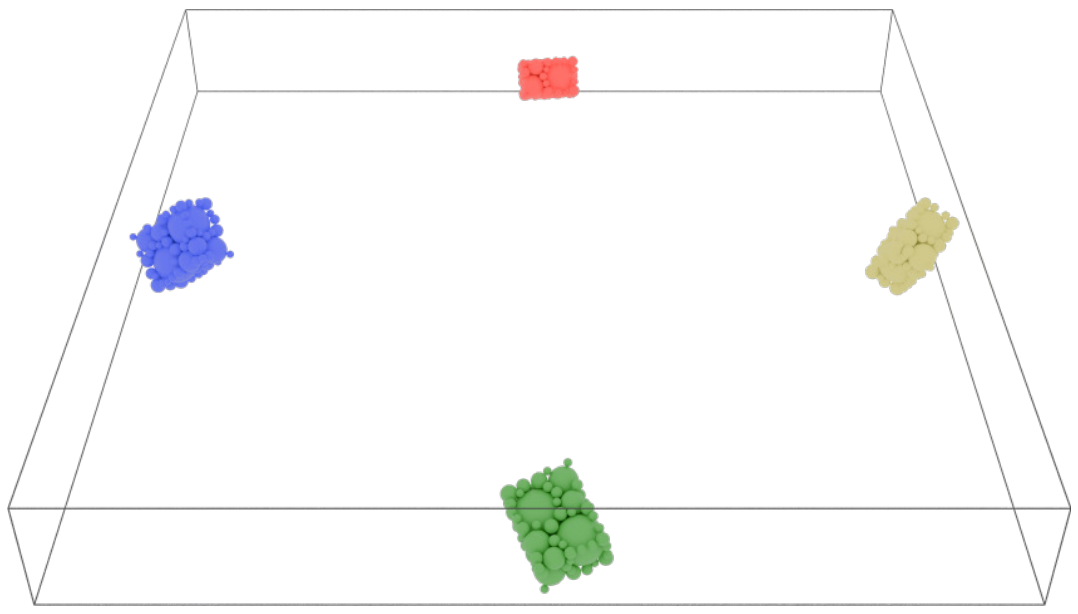

Figure 6: Initial configuration of four identical rectangular boxes in Scenario 1.

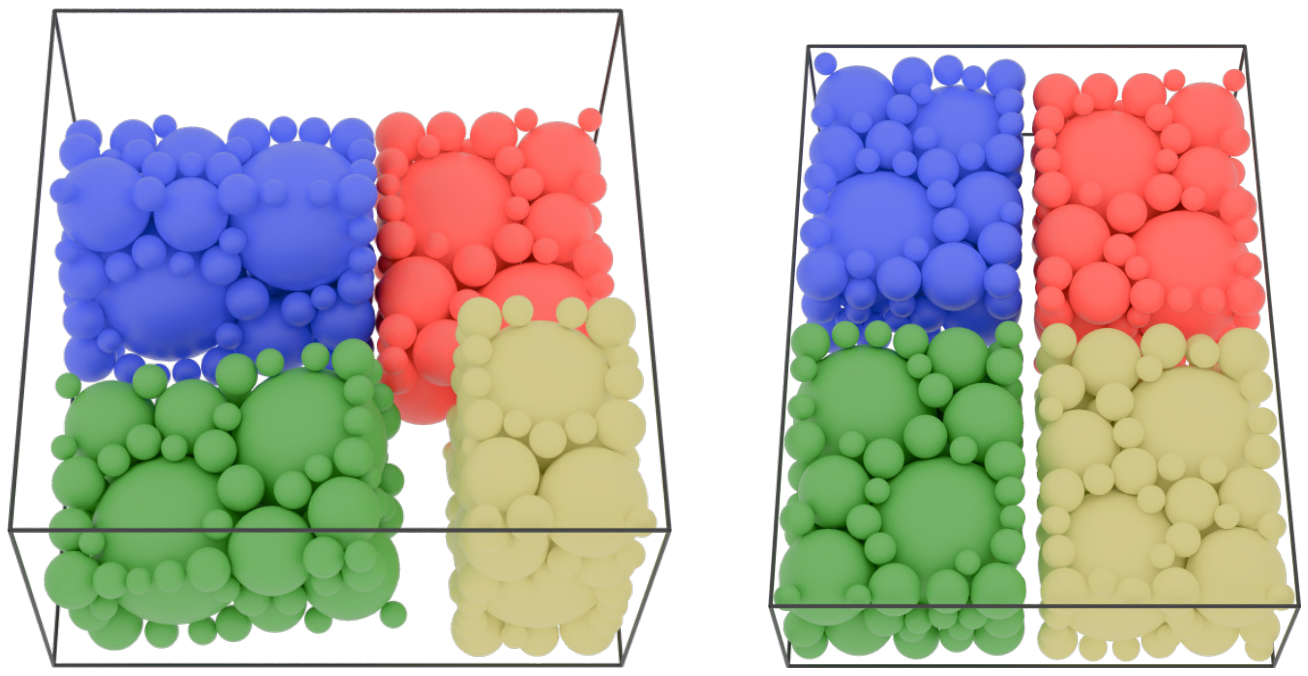

(a) GD (b) SPGD

Figure 7: Comparative final configurations for Scenario 1 by GD and SPGD. The GD configuration shows typical convergence behaviors, while SPGD demonstrates a convergence to the global optimum, representing a significantly superior solution compared to traditional GD methods.

The outcomes depicted in the figures reveal that, due to the collision constraint, GD struggled to converge to the global solution and settled in a suboptimal local minimum. In contrast, SPGD successfully converged to the global optimal configuration, effectively avoiding local minima and fulfilling the collision constraints more efficiently. To further illustrate the performance dynamics over the course of the optimization, the loss convergence history for both algorithms is plotted in Figure 8. This figure shows loss values as a function of elapsed time and the number of iterations.

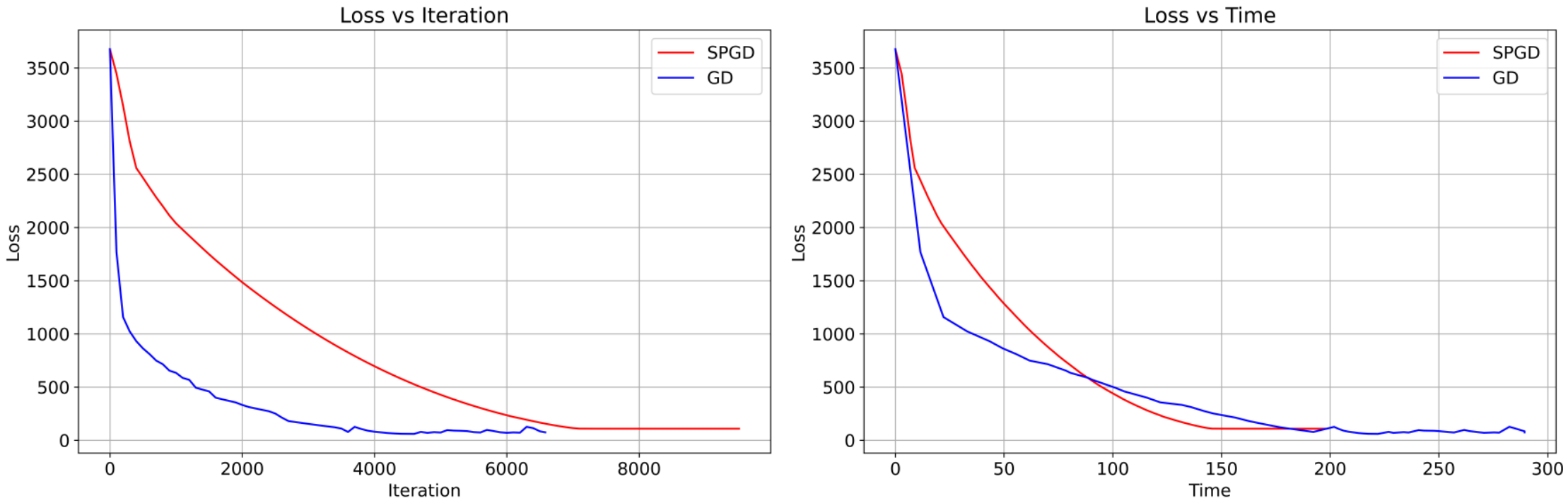

Figure 8: Loss convergence history based on elapsed time and the number of iterations for GD and SPGD in Scenario 1.

Although SPGD achieved the optimal configuration more rapidly in terms of the number of iterations, it required more computational time overall compared to GD. These plots (Figures 8) help demonstrate that while SPGD's iterations are more effective at progressing toward the global optimum, they are computationally more intensive, likely due to the complexity of the perturbation calculations and the more sophisticated collision checks involved.

This scenario underscores SPGD's strengths in effectively navigating optimization landscapes with collision constraints and its ability to reach global optima where traditional GD may fail. However, the increased computational demand highlights an area for further optimization and efficiency improvements in SPGD's implementation.

### 5.4 Analysis of Scenario 2: Eight Identical rectangular boxes

In Scenario 2, the initial configuration of eight identical rectangular boxes is depicted in Figure 9. This scenario was designed to assess each algorithm's ability to scale and manage increased numbers of objects while maintaining an efficient packing configuration. The outcomes of the final configurations found by the GD and SPGD algorithms are illustrated in Figures 10a and 10b, respectively.

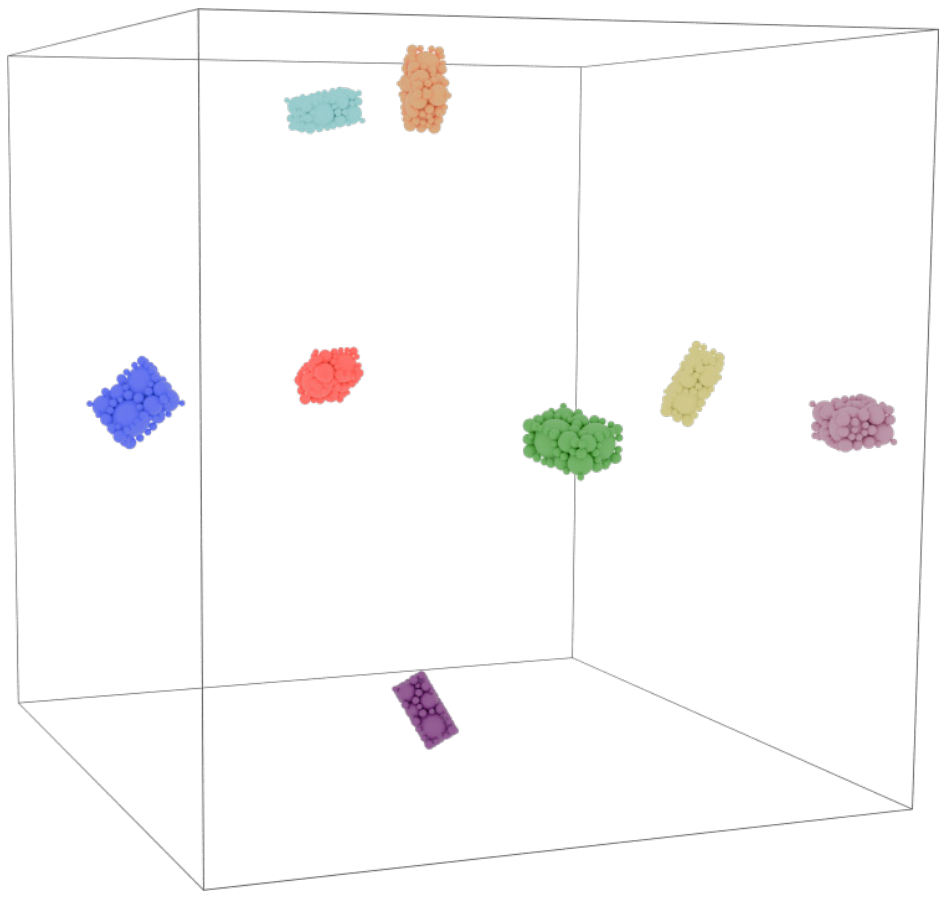

Figure 9: Initial configuration of eight identical rectangular boxes in Scenario 2.

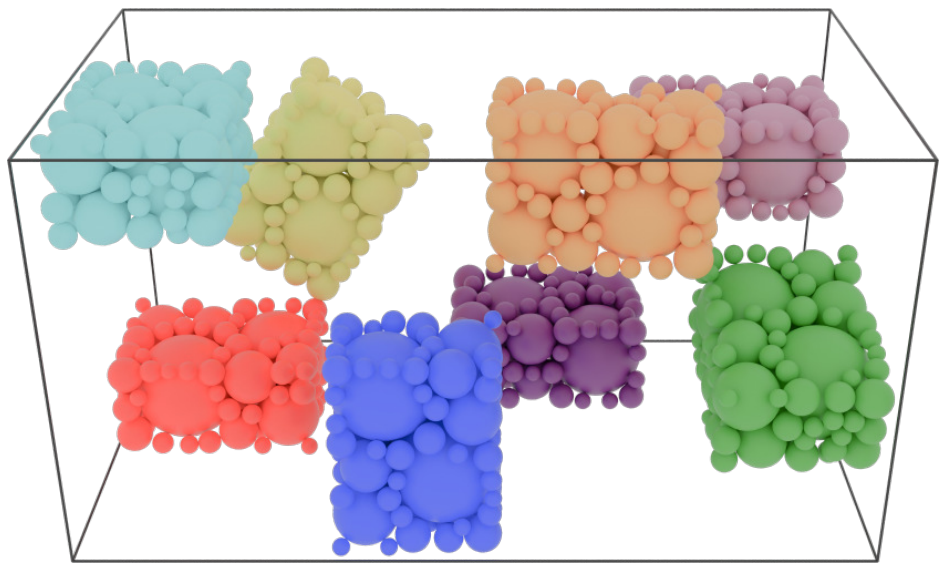

(a) GD

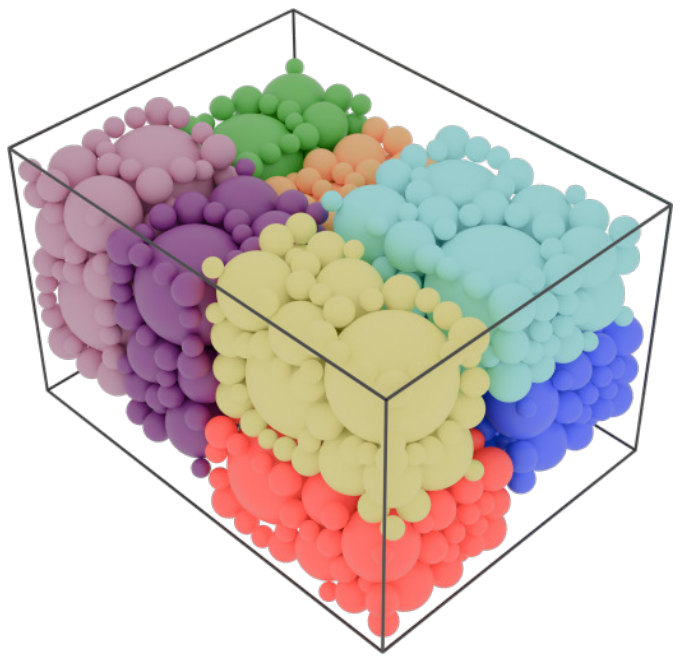

(b) SPGD

Figure 10: Comparative final configurations for Scenario 2 by GD and SPGD. GD's final arrangement demonstrates collision challenges, hindering optimal packing. In contrast, SPGD achieves a more compact configuration, effectively utilizing its adaptive perturbations to overcome collision barriers and improve packing density.

During the optimization process, the GD algorithm encountered significant issues and ceased further packing adjustments due to a collision between the yellow and red rectangular boxes, effectively stopping the optimization prematurely. In contrast, the SPGD algorithm managed to navigate around this problem and did not converge to the global optimal solution but found a notably more compact suboptimal solution, approximately three times more space-efficient than the configuration found by GD.

To further illustrate the performance dynamics over the course of the optimization, the loss convergence history for both algorithms is plotted in Figure 11, showing loss values as a function of elapsed time and the number of iterations.

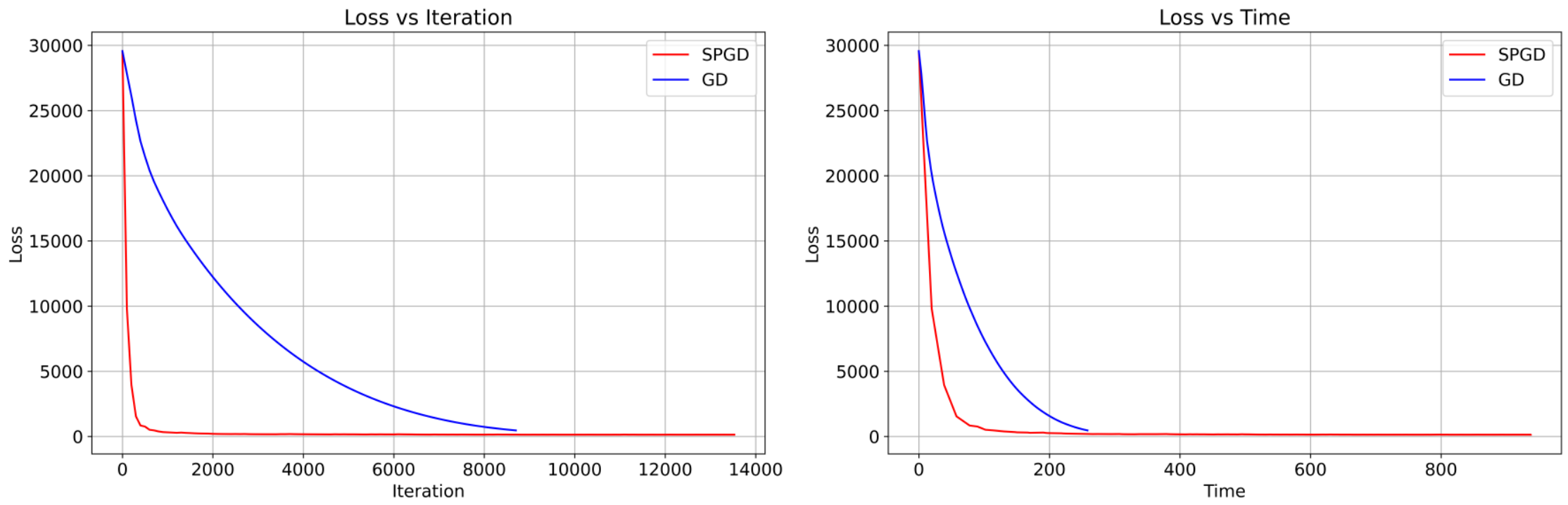

Figure 11: Loss convergence history based on elapsed time and the number of iterations for GD and SPGD in Scenario 2.

Although SPGD did not achieve the global optimum, it provided a significant improvement over GD by finding a much more compact solution rapidly. This scenario demonstrates SPGD's superior capability in effectively navigating complex landscapes and managing collision constraints dynamically compared to GD. The increased performance in finding a substantially better solution highlights the potential of SPGD for more effective space utilization in packing problems.

## 5.5 Analysis of Scenario 3: Eight rectangular boxes of Different Sizes

In Scenario 3, which introduces a higher level of complexity due to the use of eight rectangular boxes of different sizes, the initial configuration is shown in Figure 12. This setup challenges the algorithms' ability to efficiently manage and optimize space in a more heterogeneous environment.

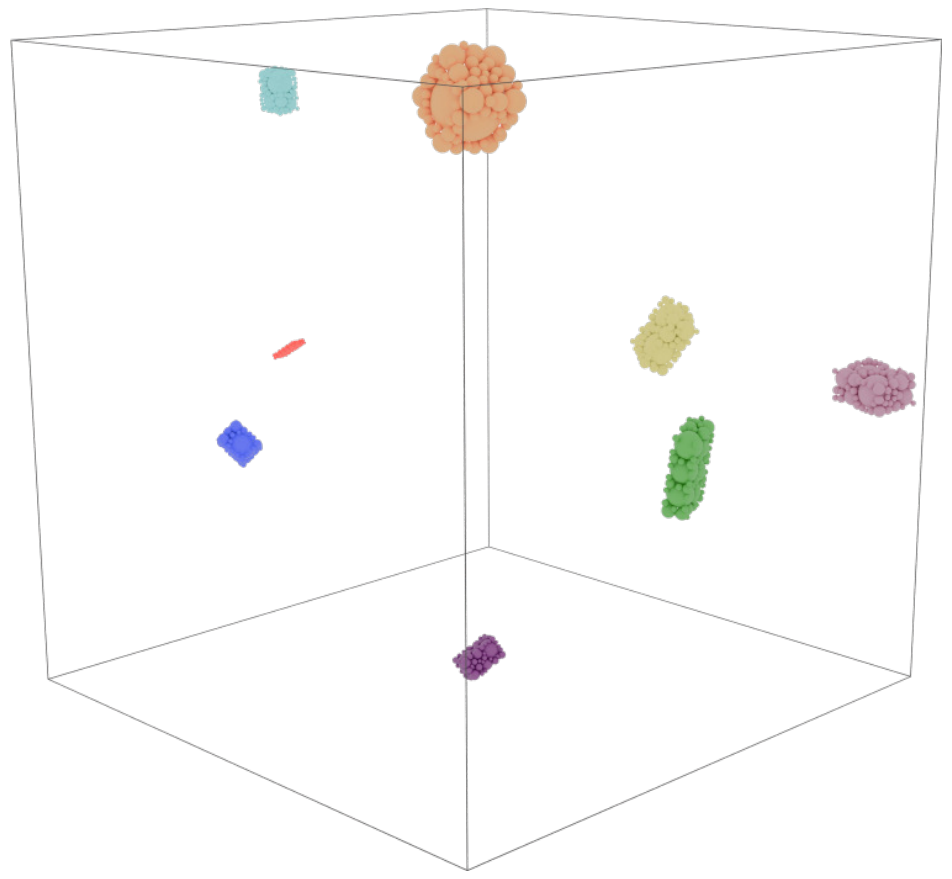

Figure 12: Initial configuration of eight rectangular boxes of different sizes in Scenario 3.

The SPGD algorithm's performance in this scenario was notably superior, as it converged to a more compact solution significantly faster than the traditional GD method. The results of the final configurations found by the GD and SPGD algorithms are shown in Figures 13a and 13b, respectively.

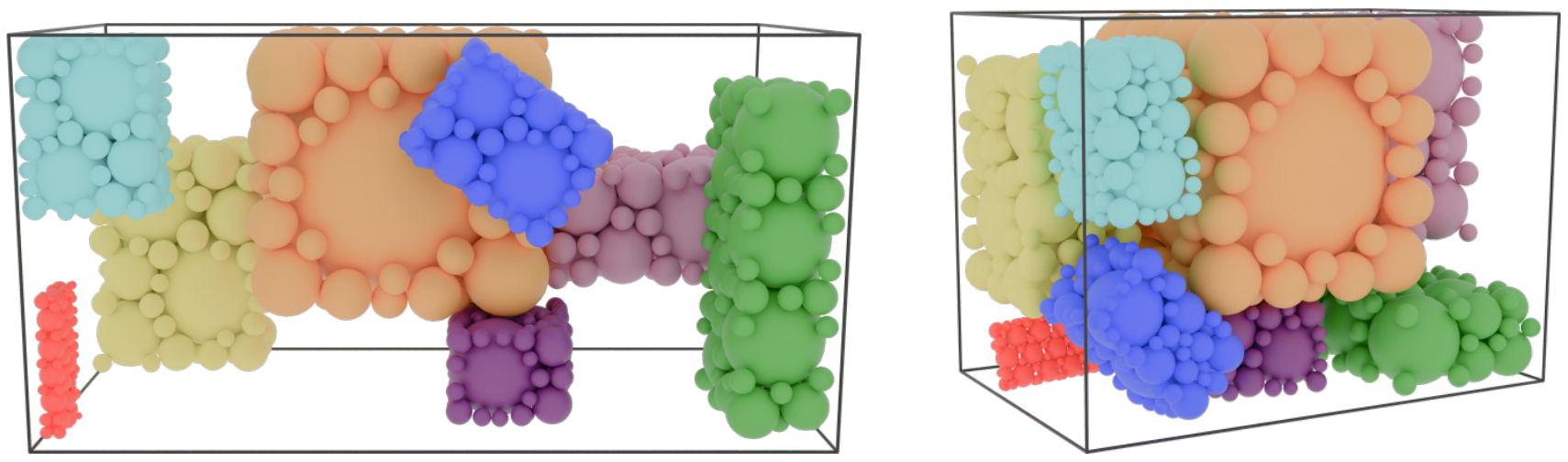

(a) GD (b) SPGD

Figure 13: Comparative final configurations for Scenario 3: Gradient Descent (left) shows less optimized packing, while Steepest Perturbed Gradient Descent (right) demonstrates a more compact and efficient arrangement.

Despite the lack of a known global optimal solution due to the varying sizes and potential configurations, SPGD effectively utilized its perturbation mechanism to explore and optimize the packing arrangement. This scenario highlights the algorithm's adaptability and efficiency in handling diverse object dimensions, which is crucial for real-world applications.

To further evaluate the performance dynamics, the loss convergence history for both algorithms is plotted in Figure 14, showing loss values as a function of elapsed time, and the number of iterations.

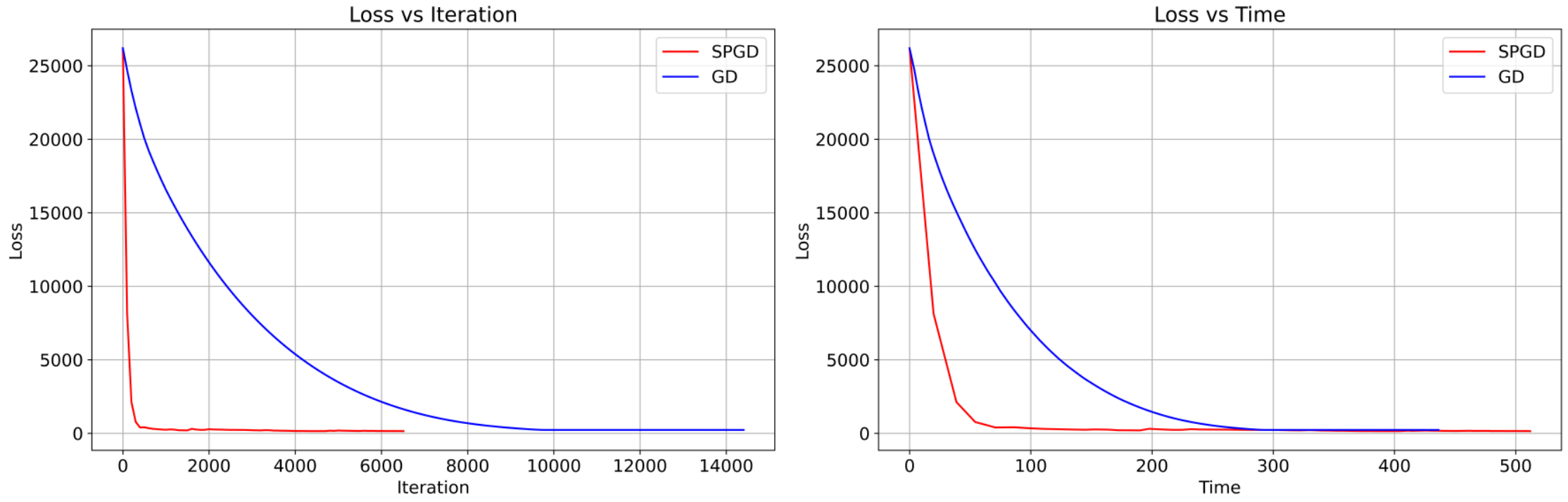

Figure 14: Loss convergence history based on elapsed time and the number of iterations for GD and SPGD in Scenario 3.

These figures demonstrate that SPGD not only achieves a more desirable outcome but also does so with greater computational efficiency in terms of iteration count, despite the complex interplay of different-sized objects. This efficiency underscores SPGD's potential as a robust tool for tackling sophisticated packing challenges where traditional methods might falter.

## 5.6 Analysis of Scenario 4: Eight Objects of Different Shapes

Scenario 4, the most complex of the scenarios tested, involved packing eight objects of different, irregular shapes such as gears, hooks, and rivets. The initial configuration is illustrated in Figure 15, which presents a diverse and challenging packing environment.

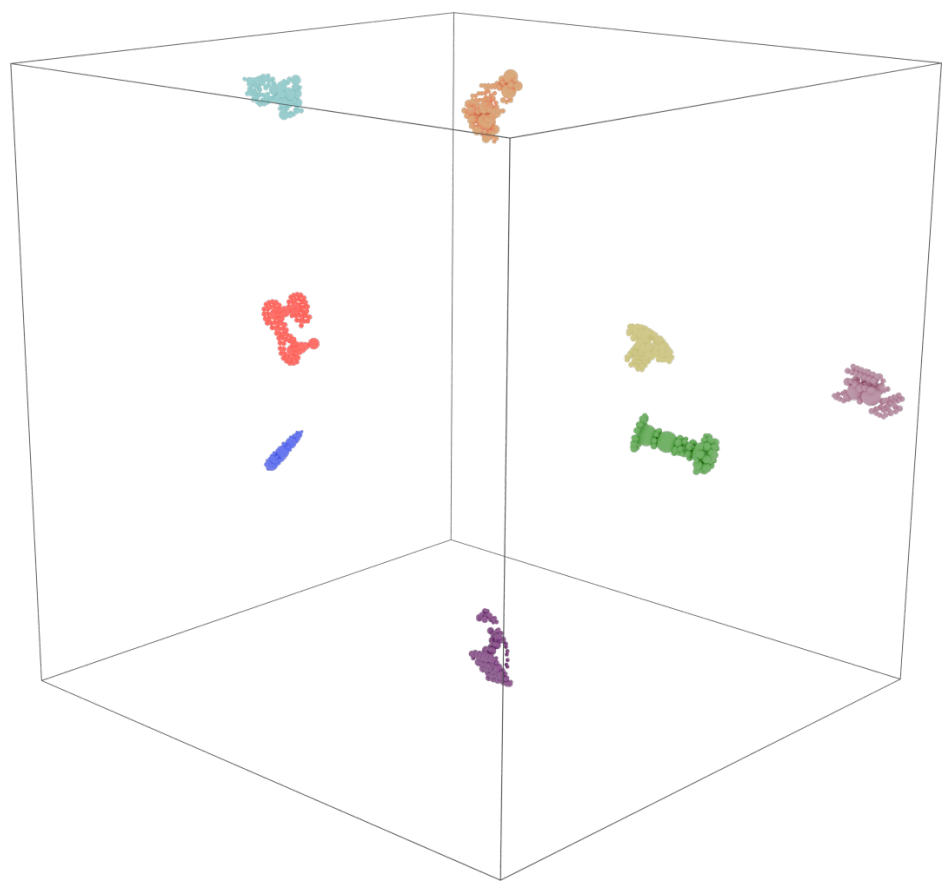

Figure 15: Initial configuration of eight objects of different shapes in Scenario 4.

In this demanding scenario, the SPGD algorithm demonstrated its robust capability by converging to a significantly more compact solution compared to the traditional GD method. Although the time taken by SPGD to find the optimal solution was comparable to that of GD, the overall optimization process required more time due to the termination condition set for no improvement in the loss value over 2000 iterations. The final configurations achieved by the GD and SPGD algorithms are shown in Figure 16a and 16b, reflecting the SPGD algorithm's effectiveness in handling complex and varied object forms.

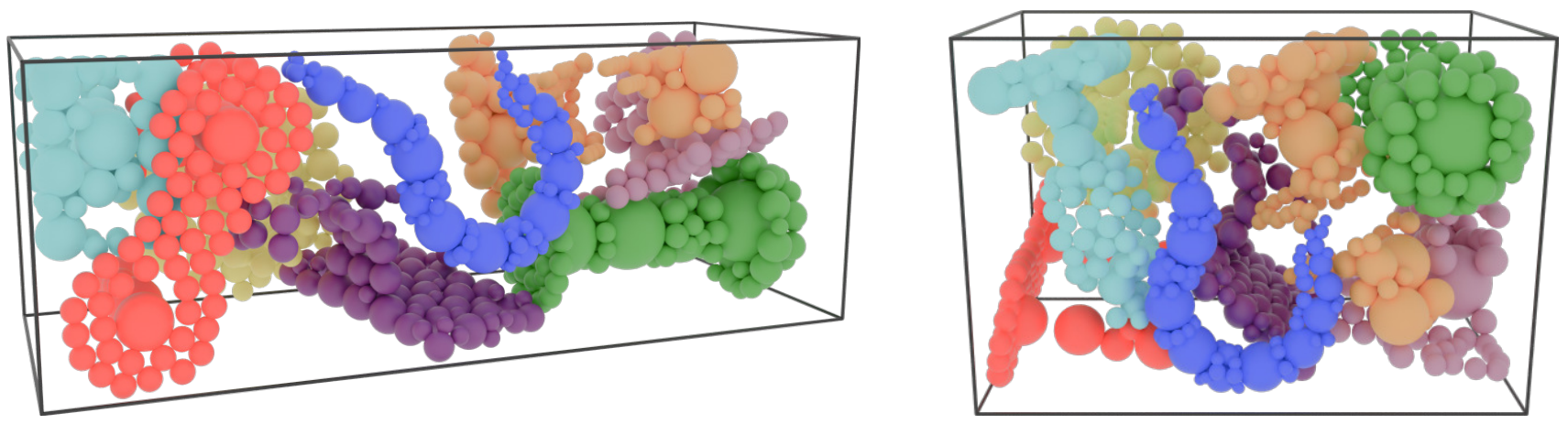

(a) GD (b) SPGD

Figure 16: Comparative final configurations for Scenario 4: Gradient Descent (a) struggles with complexity, while Steepest Perturbed Gradient Descent (b) demonstrates a significant improvement, achieving a more compact arrangement by 19.6%.

To highlight the dynamic performance of both algorithms in this scenario, Figure 17 presents the loss convergence history based on elapsed time, and the number of iterations.

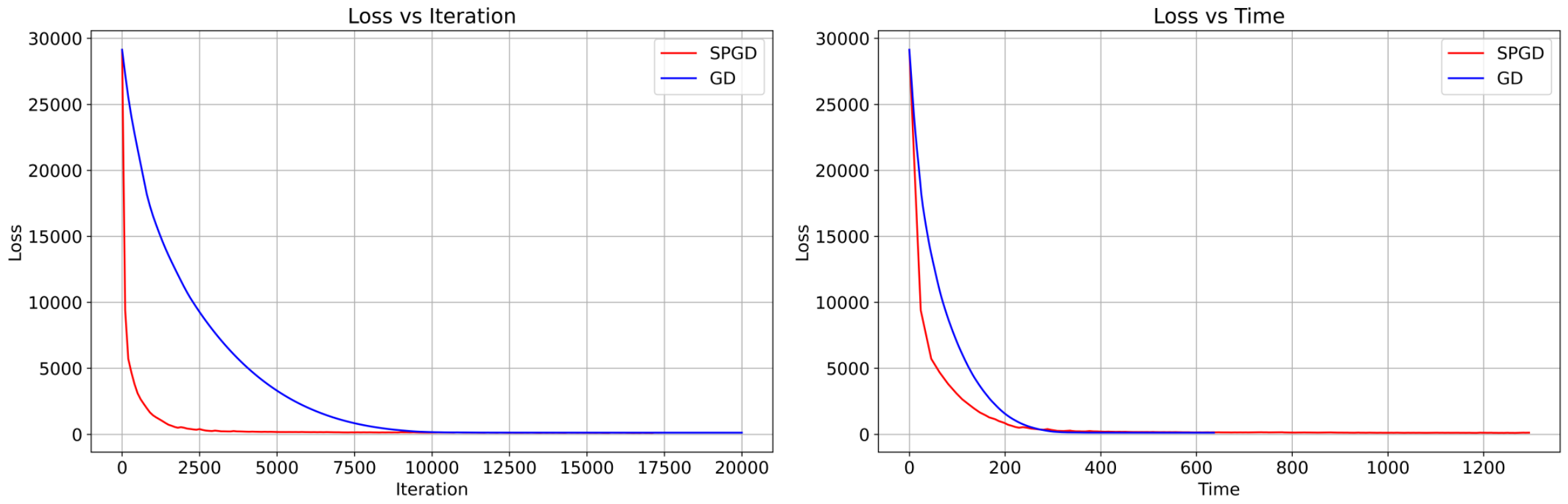

Figure 17: Loss convergence history based on elapsed time for GD and SPGD in Scenario 4.

These results underscore the SPGD algorithm's capacity to adapt to and effectively manage the intricacies of packing highly irregular objects. Although the time to reach the optimal solution was similar for both algorithms, SPGD's ability to achieve a more compact arrangement highlights its suitability for complex, real-world packing problems where shape diversity plays a critical role. The extended time required for optimization termination points to the rigorous nature of the stopping criterion, ensuring that the solution is indeed optimal before termination.

Results of these experiments are summarized in the table 9, which compares the performance of SPGD and GD in terms of best loss, and volume.

Table 9: Final Loss and Volume Comparison of SPGD and GD across different packing scenarios

| **Scenario** | **Method** | **Best Loss** | **Volume** |
|---|---|---|---|
| 1 | SPGD | 74.59 | 7.12 |
| | GD | 108.69 | 13.04 |
| 2 | SPGD | 137.38 | 16.51 |
| | GD | 465.53 | 58.24 |
| 3 | SPGD | 145.59 | 17.60 |
| | GD | 225.36 | 27.46 |
| 4 | SPGD | 103.26 | 12.34 |
| | GD | 123.40 | 14.76 |

This analysis highlights the superior adaptability and performance of SPGD, particularly in scenarios involving complex and non-uniform object configurations. The algorithm's ability to effectively shuffle and perturb component sequences contributes significantly to its success in navigating the intricate landscapes presented by these diverse packing challenges.

## 6 Conclusion

The SPGD algorithm presents a novel integration of deterministic optimization with strategic stochastic perturbations, designed to overcome the limitations of traditional gradient descent methods in non-convex landscapes and plateaus. Through comparative analyses, SPGD has demonstrated potential advantages in complex non-convex optimization challenges, consistently converging to the global optimum across 30 randomized trials per benchmark function. These results highlight both the robustness and practical utility of SPGD across a wide range of optimization scenarios.

Looking ahead, SPGD shows promise for broader applications in diverse domains and enhancements in machine learning methodologies:

- **Expanding Application Domains:** Future investigations could explore SPGD's application to fields like engineering design optimization [49], logistics, energy management, bioinformatics, and fuzzy logic parameter tuning optimization [50], showcasing its versatility and robustness.
- **Enhancements in Machine Learning:** There is potential for SPGD to significantly enhance neural network training, especially within deep learning frameworks by improving convergence rates and navigating complex parameter spaces.
- **Integration with Machine Learning Frameworks:** SPGD has already been implemented using the PyTorch framework for the 3D component packing problem, demonstrating its adaptability to complex optimization tasks. Future work could extend this integration to machine learning projects, particularly in training neural networks, thereby potentially broadening its user base and enhancing its utility in diverse applications.

- **Adaptive Perturbation Strategies:** Developing adaptive perturbation techniques that respond to specific characteristics of the optimization landscape could further refine SPGD's effectiveness, making it more problem-specific.
- **Extension to Complex Systems:** Exploring the 3D Component Packing Problem within the SPI2 framework could pave the way for handling interconnected systems with physical interactions, where topology and collision constraints add layers of complexity.

These future directions not only aim to broaden the utility of SPGD but also open new avenues for innovative research in the field of optimization.

## Acknowledgments

This work was supported in part by the National Science Foundation grants CMMI-2232612, CMMI-2312175, and by the Defense Advanced Research Projects Agency (DARPA) under the "Multi-Disciplinary Optimization for Packaging (MDOP)" program, grant number FA8750-23-C-0501.

## References

[1] Augustin Cauchy et al. Méthode générale pour la résolution des systemes d'équations simultanées. *Comp. Rend. Sci. Paris*, 25(1847):536–538, 1847.

[2] Sebastian Ruder. An overview of gradient descent optimization algorithms. *arXiv preprint arXiv:1609.04747*, 2016.

[3] Chi Jin, Rong Ge, Praneeth Netrapalli, Sham M Kakade, and Michael I Jordan. How to escape saddle points efficiently. In *International conference on machine learning*, pages 1724–1732. PMLR, 2017.

[4] Ilya Sutskever, James Martens, George Dahl, and Geoffrey Hinton. On the importance of initialization and momentum in deep learning. In *International conference on machine learning*, pages 1139–1147. PMLR, 2013.

[5] Chi Jin, Praneeth Netrapalli, and Michael I Jordan. Accelerated gradient descent escapes saddle points faster than gradient descent. In *Conference On Learning Theory*, pages 1042–1085. PMLR, 2018.

[6] Scott Kirkpatrick, C Daniel Gelatt Jr, and Mario P Vecchi. Optimization by simulated annealing. *Science*, 220(4598):671–680, 1983.

[7] John H Holland. Genetic algorithms. *Scientific American*, 267(1):66–73, 1992.

[8] Daniel Delahaye, Supatcha Chaimatanan, and Marcel Mongeau. Simulated annealing: From basics to applications. *Handbook of metaheuristics*, pages 1–35, 2019.

[9] Parsa Ghannadi, Seyed Sina Kourehli, and Seyedali Mirjalili. A review of the application of the simulated annealing algorithm in structural health monitoring (1995-2021). *Frattura ed Integrità Strutturale*, 17(64):51–76, 2023.

[10] Bobak Shahriari, Kevin Swersky, Ziyu Wang, Ryan P Adams, and Nando De Freitas. Taking the human out of the loop: A review of bayesian optimization. *Proceedings of the IEEE*, 104(1):148–175, 2016.

[11] Jasper Snoek, Hugo Larochelle, and Ryan P Adams. Practical bayesian optimization of machine learning algorithms. In *Advances in neural information processing systems*, volume 25, 2012.

[12] Peter I Frazier. A tutorial on bayesian optimization. *arXiv preprint arXiv:1807.02811*, 2018.

[13] James V Burke, Frank E Curtis, Adrian S Lewis, Michael L Overton, and Lucas EA Simões. Gradient sampling methods for nonsmooth optimization. *Numerical nonsmooth optimization: State of the art algorithms*, pages 201–225, 2020.

[14] Léon Bottou, Frank E Curtis, and Jorge Nocedal. Optimization methods for large-scale machine learning. *SIAM review*, 60(2):223–311, 2018.

[15] Ka Fai Cedric Yiu, Yanqun Liu, and Kok Lay Teo. A hybrid descent method for global optimization. *Journal of global optimization*, 28:229–238, 2004.

[16] Xin Guo, Jiequn Han, Mahan Tajrobehkar, and Wenpin Tang. Escaping saddle points efficiently with occupation-time-adapted perturbations, 2022.

[17] Xin-She Yang, TO Ting, and Mehmet Karamanoglu. Random walks, lévy flights, markov chains and metaheuristic optimization. *Future Information Communication Technology and Applications: ICFICE 2013*, pages 1055–1064, 2013.

[18] Tao Sun, Dongsheng Li, and Bao Wang. Adaptive random walk gradient descent for decentralized optimization. In *International Conference on Machine Learning*, pages 20790–20809. PMLR, 2022.

[19] David Sussillo and LF Abbott. Random walk initialization for training very deep feedforward networks. *arXiv preprint arXiv:1412.6558*, 2014.

[20] John P Cunningham and Zoubin Ghahramani. Linear dimensionality reduction: Survey, insights, and generalizations. *The Journal of Machine Learning Research*, 16(1):2859–2900, 2015.

[21] Patrick Billingsley. *Probability and measure.* John Wiley & Sons, 2017.

[22] Sanjeev Arora and Boaz Barak. *Computational complexity: a modern approach.* Cambridge University Press, 2009.

[23] MATLAB. Find minimum of constrained nonlinear multivariable function - matlab fmincon. `https://www.mathworks.com/help/optim/ug/fmincon.html`. Accessed: March 25, 2024.

[24] MATLAB. Find minimum of unconstrained nonlinear multivariable function - matlab fmincon. `https://www.mathworks.com/help/releases/R2024b/optim/ug/fminunc.html`. Accessed: April 30, 2025.

[25] MATLAB. Find minimum of function using simulated annealing algorithm. `https://www.mathworks.com/help/gads/simulannealbnd.html`. Accessed: March 25, 2024.

[26] MATLAB. Find minimum of function using simulated annealing algorithm. `https://www.mathworks.com/help/stats/bayesian-optimization-algorithm.html`. Accessed: April 30, 2025.

[27] S. Surjanovic and D. Bingham. Virtual library of simulation experiments: Test functions and datasets. Retrieved March 25, 2024, from `http://www.sfu.ca/~ssurjano`.

[28] M Hazewinkel. Theory of errors. *Encyclopedia of mathematics*, 62, 2001.

[29] The MathWorks, Inc. Matlab peaks function - matlab.

[30] Wei Wei, Dong Yang, Li Li, and Yuxuan Xia. An intravascular catheter bending recognition method for interventional surgical robots. *Machines*, 10(1), 2022.

[31] David Ackley. *A connectionist machine for genetic hillclimbing*, volume 28. Springer Science & Business Media, 2012.

[32] Eric E Easom. *A survey of global optimization techniques.* PhD thesis, University of Louisville, 1990.

[33] A. V. Levy and Antonio Montalvo. The tunneling algorithm for the global minimization of functions. *Siam Journal on Scientific and Statistical Computing*, 6:15–29, 1985.

[34] Anne Gershenson, Shachi Gosavi, Pietro Faccioli, and Patrick L Wintrode. Successes and challenges in simulating the folding of large proteins. *Journal of Biological Chemistry*, 295(1):15–33, 2020.

[35] Ken Dill and Hue Chan. From Levinthal to pathways to funnels. *Nature structural biology*, 4:10–9, 02 1997.

[36] Zahra Shahbazi, Horea T. Ilieş, and Kazem Kazerounian. Hydrogen Bonds and Kinematic Mobility of Protein Molecules. *Journal of Mechanisms and Robotics*, 2(2):021009, 04 2010.

[37] Christopher Madden, Peter Bohnenkamp, Kazem Kazerounian, and Horea T Ilieş. Residue level three-dimensional workspace maps for conformational trajectory planning of proteins. *The International Journal of Robotics Research*, 28(4):450–463, 2009.

[38] Pouya Tavousi, Morad Behandish, Horea T Ilieş, and Kazem Kazerounian. Protofold ii: Enhanced model and implementation for kinetostatic protein folding. *Journal of Nanotechnology in Engineering and Medicine*, 6(3):034601, 2015.

[39] Alireza Mohammadi and Mohammad Al Janaideh. Sign gradient descent algorithms for kinetostatic protein folding. In *2023 International Conference on Manipulation, Automation and Robotics at Small Scales (MARSS)*, pages 1–6. IEEE, 2023.

[40] Satya RT Peddada, Lawrence E Zeidner, Horea T Ilies, Kai A James, and James T Allison. Toward holistic design of spatial packaging of interconnected systems with physical interactions (spi2). *Journal of Mechanical Design*, 144(12):120801, 2022.

[41] S. Szykman and J. Cagan. A Simulated Annealing-Based Approach to Three-Dimensional Component Packing. *Journal of Mechanical Design*, 117(2A):308–314, 06 1995.

[42] Satya R. T. Peddada, Lawrence E. Zeidner, Horea T. Ilies, Kai A. James, and James T. Allison. Toward Holistic Design of Spatial Packaging of Interconnected Systems With Physical Interactions (SPI2). *Journal of Mechanical Design*, 144(12):120801, 08 2022.

[43] Mohammad M. Behzadi, Peter Zaffetti, Jiangce Chen, Lawrence E. Zeidner, and Horea T Ilieş. Spatial component packing and routing optimization with physical interaction using maximal disjoint ball decomposition. *Journal of Mechanical Design*, 2024. in press.

[44] Mikola Lysenko. Fourier collision detection. *The International Journal of Robotics Research*, 32(4):483–503, 2013.

[45] Morad Behandish and Horea T. Ilieş. Peg-in-Hole Revisited: A Generic Force Model for Haptic Assembly. *Journal of Computing and Information Science in Engineering*, 15(4):041004, 08 2015.

[46] Morad Behandish and Horea T Ilieş. Analytic methods for geometric modeling via spherical decomposition. *Computer-Aided Design*, 70:100–115, 2016.

[47] Lydia E Kavraki. Computation of configuration-space obstacles using the fast fourier transform. *IEEE Transactions on Robotics and Automation*, 11(3):408–413, 1995.

[48] Qiaodong Cui, Victor Rong, Desai Chen, and Wojciech Matusik. Dense, interlocking-free and scalable spectral packing of generic 3D objects. *ACM Trans. Graph.*, 42(4):141–1, 2023.

[49] Shanglong Zhang and Julián A Norato. Finding better local optima in topology optimization via tunneling. In *International Design Engineering Technical Conferences and Computers and Information in Engineering Conference*, volume 51760, page V02BT03A014. American Society of Mechanical Engineers, 2018.

[50] Amir Mohammad Vahedi, Hadi Nobahari, and Meysam Alizad. Fuzzy gain scheduling of artificial potential fields for online path planning and obstacle avoidance of an aerial robot. In *2022 10th RSI International Conference on Robotics and Mechatronics (ICRoM)*, pages 309–316, 2022.

# A Appendix: Descent Lemma

**Lemma A.1** (Descent lemma for $\ell$-smooth $f$)**.** *If $f : \mathbb{R}^d \to \mathbb{R}$ has $\ell$-Lipschitz gradient, then for all $x, y \in \mathbb{R}^d$,*

$$f(\mathbf{y}) \;\le\; f(\mathbf{x}) \;+\; \nabla f(\mathbf{x})^\top(\mathbf{y}-\mathbf{x}) \;+\; \frac{\ell}{2}\,\|\mathbf{y}-\mathbf{x}\|_2^2. \tag{10}$$

This inequality can be viewed as a quadratic approximation to $f$ around $\mathbf{x}$.

**Corollary A.2** (One-step decrease of GD (Lemma 9 in [3]))**.** *Let $\mathbf{x}^+ = \mathbf{x} - \alpha\,\nabla f(\mathbf{x})$ with step size $\alpha > 0$. Under the assumptions of A.1,*

$$f(\mathbf{x}^+) \;\le\; f(\mathbf{x}) \;-\; \alpha\Big(1 - \tfrac{\alpha\ell}{2}\Big)\,\|\nabla f(\mathbf{x})\|_2^2. \tag{11}$$

*In particular, if $\alpha \le 1/\ell$ then*

$$f(\mathbf{x}^+) \;\le\; f(\mathbf{x}) \;-\; \frac{\alpha}{2}\,\|\nabla f(\mathbf{x})\|_2^2. \tag{12}$$

*If $\alpha = c/\ell$ with a universal $c \in (0, 1]$, then*

$$f(\mathbf{x}^+) \;\le\; f(\mathbf{x}) \;-\; \frac{c}{2\ell}\,\|\nabla f(\mathbf{x})\|_2^2. \tag{13}$$

**Proof A.3** (Proof of A.2)**.** *Apply* (10) *with* $\mathbf{y} = \mathbf{x} - \alpha\nabla f(\mathbf{x})$*:*

$$f(\mathbf{x}^+) \le f(\mathbf{x}) - \alpha\|\nabla f(\mathbf{x})\|_2^2 + \frac{\ell}{2}\alpha^2\|\nabla f(\mathbf{x})\|_2^2 = f(\mathbf{x}) - \alpha\Big(1 - \tfrac{\alpha\ell}{2}\Big)\|\nabla f(\mathbf{x})\|_2^2,$$

*which gives* (11)*;* (12) *follows since $\alpha\ell \le 1$, and* (13) *follows by substituting $\alpha = c/\ell$.*

# B Appendix: Notation and Symbols

**Notation**

- $d \in \mathbb{N}$ is the space dimension; vectors are in $\mathbb{R}^d$. $\|\cdot\|_2$: Euclidean norm.
- For $A \in \mathbb{R}^{d\times d}$, $\|A\|_{\mathrm{op}}$ is the spectral norm; $\lambda_{\min}(A)$: smallest eigenvalue.
- $\mathbb{B}(0,r) = \{u \in \mathbb{R}^d : \|u\|_2 \le r\}$; $\mathrm{Unif}(\mathbb{B}(0,r))$: uniform distribution within a ball of radius $r$.
- $\nabla f(\mathbf{x})$, $\nabla^2 f(\mathbf{x})$: gradient and Hessian of $f$ at $\mathbf{x}$.
- $\tilde{O}(\cdot)$: big-$O$ up to polylogarithms in natural parameters (e.g., $d$, $1/\delta$).
- An *iteration* is one GD step; a *perturbation round* is when seeds are injected and one is selected.

| **Symbol** | Meaning |
|---|---|
| $d$ | Dimension. |
| $\ell$ | Gradient Lipschitz constant. |
| $\rho$ | Hessian Lipschitz constant. |
| $f^\star$ | Infimum of $f$. |
| $\Delta_f$ | $f(x_0) - f^\star$ (initial suboptimality). |
| $\alpha$ | Step size, $= c/\ell$. |
| $\epsilon$ | Gradient tolerance in SOSP. |
| $\sqrt{\rho\epsilon}$ | Curvature tolerance in SOSP. |
| Amp | Perturbation radius. |
| $\tau$ | Escape horizon (and period in rollout analysis). |
| $\mathrm{Iter}_P$ | Period between perturbations. |
| $N_P$ | Number of seeds per round. |
| $p_0$ | Per-seed success probability (PGD), $\ge 1/\mathrm{poly}(d)$. |
| $\delta$ | Overall failure budget. |
| $\delta_{\mathrm{round}}$ | Per-round failure budget ($= \delta/R$). |
| $T$ | Total GD iterations. |
| $R$ | Number of rounds ($\asymp T/\mathrm{Iter}_P$). |